\documentclass[12pt]{article} 
\usepackage{latexsym}
\usepackage{amsmath}
\usepackage{bm}
\usepackage{graphicx,color}
\usepackage{amssymb}
\usepackage{cite}
\usepackage{amsfonts, amscd, amsthm}

\newcommand{\sn}{{\rm sn}}
\newcommand{\cn}{{\rm cn}}
\newcommand{\dn}{{\rm dn}}

\begin{document} 

\vspace{15mm}  

\centerline{\Large \bf  Full Addition Formulae of Genus Two Hyperelliptic Functions}
\centerline{\Large \bf  by the Duplication Method}
\vspace{10mm}

\centerline {{\it Tezukayama University , Tezukayama 7, Nara 631, Japan
}}

\centerline{Kazuyasu Shigemoto\footnote{E-mail address:
shigemot@tezukayama-u.ac.jp}}

\vspace{10mm}

\centerline{\bf Abstract}  
\noindent 
We have given full algebraic addition formulae of genus two 
hyperelliptic functions by the duplication method.
This full addition formulae
according to the duplication method give some hints of Lie group structure
in addition formulae of genus two hyperelliptic functions.

\vspace{5mm}

\noindent
\vspace{10mm}

\section{Introduction}
\setcounter{equation}{0}
Addition formulae of genus two hyperelliptic theta functions 
were first obtained by Rosenhain\cite{Rosenhain1,Rosenhain2,Konigsberger}  
by using the Riemann's theta identity, that is, the identity of the sum of the  
$4$-product of hyperelliptic theta functions. 
Full addition formulae of genus two hyperelliptic functions were
given by Kossak\cite{Kossak} but he gave only the 
sketch to derive full addition formulae but did not give full
addition formulae. In our previous paper\cite{Shigemoto4} , 
we gave full addition formulae of genus two hyperelliptic functions
according to Kossak's scenario.\\
While G\"{o}pel \cite{Gopel1,Gopel2}  obtained addition 
formulae of hyperelliptic 
theta functions by using the duplication method, that is, the identity
 of the sum of the $2$-product of the hyperelliptic theta functions.
Then by using G\"{o}pel method, we will give another expression of 
full addition formulae of genus two hyperelliptic functions.\\

\vspace{5mm}
\setcounter{equation}{0}
\section{Addition Formulae of Genus Two Hyperelliptic Functions}
Here will give full addition formulae of genus two hyperelliptic functions
by using the duplication method.\\
%
\subsection{Duplication Formulae of Genus Two Hyperelliptic 
Theta Functions}\
Genus two hyperelliptic theta functions with two variables are defined by  
\begin{eqnarray}
&&\hskip -15mm  \vartheta[\begin{array}{cc} a \ c \\ b \ d \\ 
\end{array}](x,y;\tau_1, \tau_2, \tau_{12})
\nonumber\\
&&\hskip -15mm =\sum_{m,n\in Z} \exp 
\left\{ \pi i  
\Bigl( \tau_1 (m+\frac{a}{2})^2+\tau_2 (n+\frac{c}{2})^2+
2 \tau_{12} (m+\frac{a}{2}) (n+\frac{c}{2}) \Bigr) \right.
\nonumber\\
&&\hskip -15mm \left. +2\pi i \ 
\Bigl( (m+\frac{a}{2})(x+\frac{b}{2})+(n+\frac{c}{2})(y+\frac{d}{2})
\Bigr)
\right\}  ,
\label{2e1}
\end{eqnarray}
where we assume that ${\rm Im} \tau_1>0$, ${\rm Im}\tau_2>0$,
$({\rm Im} \tau_1)({\rm Im}\tau_2)-({\rm Im} \tau_{12})^2>0$ in order
that the summation of $m,n \in Z$ becomes convergent. 
We can rename $m\rightarrow m$, $n\rightarrow -n$, so that we can always choose
${\rm Im} \tau_{12}>0$.\\
We simply write
 $\displaystyle{ \vartheta[
\begin{array}{cc} a \ c \\ b \ d \\ \end{array}](x,y)
=\vartheta[\begin{array}{cc} a \ c \\ b \ d \\ 
\end{array}](x,y;\tau_1, \tau_2, \tau_{12})}$,
which has the property 
\begin{eqnarray}
&&\vartheta[\begin{array}{cc} a \ c \\ b \ d \\ \end{array}](x,y)
=\vartheta[\begin{array}{cc} a+2 &\hspace{-2mm} c \\ b &\hspace{-2mm} d \\ \end{array}](x,y)
=\vartheta[\begin{array}{cc} a &\hspace{-2mm} c+2 \\ b &\hspace{-2mm} d \\ \end{array}](x,y) 
=(-1)^a \vartheta[\begin{array}{cc} a &\hspace{-2mm} c \\ b+2 &\hspace{-2mm} d \\ \end{array}](x,y)\nonumber\\
&&=(-1)^c \vartheta[\begin{array}{cc} a &\hspace{-2mm} c \\ b &\hspace{-2mm} d+2 \\ \end{array}](x,y)   .
\nonumber
\end{eqnarray}
We also use the notation 
 $\displaystyle{ \varTheta[\begin{array}{cc} a \ c \\ b \ d \\ \end{array}](x,y)
=\vartheta[\begin{array}{cc} a \ c \\ b \ d \\ 
\end{array}](x,y;2\tau_1, 2\tau_2, 2\tau_{12})}$.\\ 
Original duplication formulae are given by G\"{o}pel\cite{Gopel1,Gopel2,Borchardt,Cayley}, 
but we use generalized duplication formulae of the form  
\begin{eqnarray}
&&\vartheta[\begin{array}{cc} a \ c \\ b \ d \\ \end{array}](x_1,y_1) \ 
\vartheta[\begin{array}{cc} e \ g \\ f \ h \\ \end{array}](x_2,y_2) 
\nonumber\\
&&
=\varTheta[\begin{array}{cc} \frac{a+e}{2} & \frac{c+g}{2} \\ b+f & d+h \\ \end{array}](x_1+x_2 , y_1+y_2) \ 
\varTheta[\begin{array}{cc} \frac{a-e}{2} & \frac{c-g}{2} \\ b-f & d-h \\ \end{array}](x_1-x_2 , y_1-y_2) 
\nonumber\\
&&+\varTheta[\begin{array}{cc} \frac{a+e}{2} & \frac{c+g+2}{2} \\ b+f & d+h \\ \end{array}](x_1+x_2 , y_1+y_2) \  
\varTheta[\begin{array}{cc} \frac{a-e}{2} & \frac{c-g+2}{2} \\ b-f & d-h \\ \end{array}](x_1-x_2 , y_1-y_2) \
\nonumber\\
&&+\varTheta[\begin{array}{cc} \frac{a+e+2}{2} & \frac{c+g}{2} \\ b+f & d+h \\ \end{array}](x_1+x_2 , y_1+y_2) \ 
\varTheta[\begin{array}{cc} \frac{a-e+2}{2} & \frac{c-g}{2} \\ b-f & d-h \\ \end{array}](x_1-x_2 , y_1-y_2) 
\nonumber\\
&&+\varTheta[\begin{array}{cc} \frac{a+e+2}{2} & \frac{c+g+2}{2} \\ b+f & d+h \\ 
\end{array}](x_1+x_2 , y_1+y_2) \ 
\varTheta[\begin{array}{cc} \frac{a-e+2}{2} & \frac{c-g+2}{2} \\ b-f & d-h \\ \end{array}](x_1-x_2 , y_1-y_2) .
\label{2e2}
\end{eqnarray}
In the product of $\vartheta[\begin{array}{cc} * \ * \\ * \ * \\ \end{array}](x_1,y_1)$ 
and $\vartheta[\begin{array}{cc} * \ * \\ * \ * \\ \end{array}](x_2,y_2)$, the quadratic form 
of the exponent is decomposed into the symmetric part $(x_1+x_2,y_1+y_2)$ and the 
anti-symmetric part $(x_1-x_2,y_1-y_2)$ with various characteristics.
Putting $x_1=u_1+u_2$, $x_2=u_1-u_2$, $y_1=v_1+v_2$, $y_2=v_1-v_2$, we rewrite above formulae
of the form
\begin{eqnarray}
&&\vartheta[\begin{array}{cc} a \ c \\ b \ d \\ \end{array}](u_1+u_2 ,v_1+v_2) \ 
\vartheta[\begin{array}{cc} e \ g \\ f \ h \\ \end{array}](u_1-u_2, v_1-v_2) 
\nonumber\\
&&
=\varTheta[\begin{array}{cc} \frac{a+e}{2} & \frac{c+g}{2} \\ b+f & d+h \\ \end{array}](2u_1 , 2v_1) \ 
\varTheta[\begin{array}{cc} \frac{a-e}{2} & \frac{c-g}{2} \\ b-f & d-h \\ \end{array}](2u_2 , 2v_2) 
\nonumber\\
&&+\varTheta[\begin{array}{cc} \frac{a+e}{2} & \frac{c+g+2}{2} \\ b+f & d+h \\ \end{array}](2u_1 , 2v_1) \
\varTheta[\begin{array}{cc} \frac{a-e}{2} & \frac{c-g+2}{2} \\ b-f & d-h \\ \end{array}](2u_2 , 2v_2) \
\nonumber\\
&&+\varTheta[\begin{array}{cc} \frac{a+e+2}{2} & \frac{c+g}{2} \\ b+f & d+h \\ \end{array}](2u_1 , 2v_1) \
\varTheta[\begin{array}{cc} \frac{a-e+2}{2} & \frac{c-g}{2} \\ b-f & d-h \\ \end{array}](2u_2 , 2v_2) 
\nonumber\\
&&+\varTheta[\begin{array}{cc} \frac{a+e+2}{2} & \frac{c+g+2}{2} \\ b+f & d+h \\ \end{array}](2u_1 , 2v_1) \ 
\varTheta[\begin{array}{cc} \frac{a-e+2}{2} & \frac{c-g+2}{2} \\ b-f & d-h \\ \end{array}](2u_2 , 2v_2) .
\label{2e3}
\end{eqnarray}\\
In order to obtain addition formulae of hyperelliptic functions, 
we put $e=a$, $g=c$, $f=0$, $h=0$ in Eq.(\ref{2e3}) and we obtain
\begin{eqnarray}
&&\vartheta[\begin{array}{cc} a \ c \\ b \ d \\ \end{array}](u_1+u_2 ,v_1+v_2) \ 
\vartheta[\begin{array}{cc} a \ c \\ 0 \ 0 \\ \end{array}](u_1-u_2, v_1-v_2) 
\nonumber\\
&&
=\varTheta[\begin{array}{cc} a \ c \\ b \ d \\ \end{array}](2u_1 , 2v_1) \  
\varTheta[\begin{array}{cc} 0 \ 0 \\ b \ d \\ \end{array}](2u_2 , 2v_2) 
\nonumber\\
&&+\varTheta[\begin{array}{cc} a &\hspace{-2mm} c+1 \\  b &\hspace{-2mm}  d \\ \end{array}](2u_1 , 2v_1) \
\varTheta[\begin{array}{cc} 0 \ 1 \\ b \ d \\ \end{array}](2u_2 , 2v_2) \
\nonumber\\
&&+\varTheta[\begin{array}{cc} a+1 &\hspace{-2mm}  c \\ b &\hspace{-2mm}  d \\ \end{array}](2u_1 , 2v_1) \ 
\varTheta[\begin{array}{cc} 1 \ 0 \\ b \ d \\ \end{array}](2u_2 , 2v_2) 
\nonumber\\
&&+\varTheta[\begin{array}{cc} a+1 &\hspace{-2mm}  c+1 \\ b &\hspace{-2mm}  d \\ \end{array}](2u_1 , 2v_1) \ 
\varTheta[\begin{array}{cc} 1 \ 1 \\ b \ d \\ \end{array}](2u_2 , 2v_2) .
\label{2e4}
\end{eqnarray}\\
In order to obtain functional relations(I), we put $e=a$, $f=b$, $g=c$, $h=d$, $u_2=0$, $v_2=0$, $u_1=u$, $v_1=v$ in Eq.(\ref{2e3}) and we obtain 
\begin{eqnarray}
&&\hskip-15mm \vartheta[\begin{array}{cc} a \ c \\ b \ d \\ \end{array}]^2(u ,v)
\nonumber\\
&&\hskip-15mm 
=\varTheta[\begin{array}{cc} a \ c \\ 2b \ 2d \\ \end{array}](2u , 2v) \ 
\varTheta[\begin{array}{cc} 0 \ 0 \\ 0 \ 0 \\ \end{array}](0 , 0) 
+\varTheta[\begin{array}{cc} a &\hspace{-2mm}  c+1 \\ 2b &\hspace{-2mm}  2d \\ \end{array}](2u , 2v) \
\varTheta[\begin{array}{cc} 0 \ 1 \\ 0 \ 0 \\ \end{array}](0 , 0) \
\nonumber\\
&&\hskip-15mm 
+\varTheta[\begin{array}{cc} a+1 &\hspace{-2mm}  c \\ 2b &\hspace{-2mm}  2d \\ \end{array}](2u , 2v) \
\varTheta[\begin{array}{cc} 1 \ 0 \\ 0 \ 0 \\ \end{array}](0 , 0) 
+\varTheta[\begin{array}{cc} a+1 &\hspace{-2mm}  c+1 \\ 2b &\hspace{-2mm}  2d \\ \end{array}](2u , 2v) \
\varTheta[\begin{array}{cc} 1 \ 1 \\ 0 \ 0 \\ \end{array}](0 , 0)   .
\label{2e5}
\end{eqnarray}
This is the original G\"{o}pel's duplication formula.
\\
In order to obtain functional relations(II), we put $e=a$, $g=c$, $f=0$, $h=0$, $u_2=0$, $v_2=0$, $u_1=u$, $v_1=v$ and we obtain
\begin{eqnarray}
&&\hskip-15mm \vartheta[\begin{array}{cc} a \ c \\ b \ d \\ \end{array}](u , v) \ 
\vartheta[\begin{array}{cc} a \ c \\ 0 \ 0 \\ \end{array}](u, v) 
\nonumber\\
&&\hskip-15mm
=\varTheta[\begin{array}{cc} a \ c \\ b \ d \\ \end{array}](2u , 2v) \  
\varTheta[\begin{array}{cc} 0 \ 0 \\ b \ d \\ \end{array}](0 , 0) 
+\varTheta[\begin{array}{cc} a &\hspace{-2mm} c+1 \\  b &\hspace{-2mm}  d \\ \end{array}](2u , 2v) \
\varTheta[\begin{array}{cc} 0 \ 1 \\ b \ d \\ \end{array}](0 , 0) \
\nonumber\\
&&\hskip-15mm 
+\varTheta[\begin{array}{cc} a+1 &\hspace{-2mm}  c \\ b &\hspace{-2mm}  d \\ \end{array}](2u , 2v) \ 
\varTheta[\begin{array}{cc} 1 \ 0 \\ b \ d \\ \end{array}](0 , 0) 
+\varTheta[\begin{array}{cc} a+1 &\hspace{-2mm}  c+1 \\ b &\hspace{-2mm}  d \\ \end{array}](2u , 2v) \ 
\varTheta[\begin{array}{cc} 1 \ 1 \\ b \ d \\ \end{array}](0 , 0)   .
\label{2e6}
\end{eqnarray}

\subsection{Twelve Addition Formulae}
For genus two case, there are sixteen hyperelliptic theta functions, so that there are 
fifteen hyperelliptic functions. Then there are fifteen addition formulae of hyperelliptic 
functions. Twelve addition formulae of full fifteen addition formulae are easily obtained 
by this duplication method by using Eq.(\ref{2e4}).\\
\\
\noindent
{\bf A) sector:\ $a=0$, $c=0$;}\\
In this sector of Eq.(\ref{2e4}), addition formulae for
$\vartheta[\begin{array}{cc} 0 \ 0 \\ * \ * \\ \end{array}](u_1+u_2 ,v_1+v_2)$
are obtained by using following formulae
\begin{eqnarray}
&&\hskip-15mm 
\vartheta[\begin{array}{cc} 0 \ 0 \\ b \ d \\ \end{array}](u_1+u_2 ,v_1+v_2) \ 
\vartheta[\begin{array}{cc} 0 \ 0 \\ 0 \ 0 \\ \end{array}](u_1-u_2, v_1-v_2) 
\nonumber\\
&&\hskip-15mm 
=\varTheta[\begin{array}{cc} 0 \ 0 \\ b \ d \\ \end{array}](2u_1 , 2v_1) \  
\varTheta[\begin{array}{cc} 0  \ 0 \\ b \ d \\ \end{array}](2u_2 , 2v_2) 
+\varTheta[\begin{array}{cc} 0 \ 1 \\ b \ d \\ \end{array}](2u_1 , 2v_1) \ 
\varTheta[\begin{array}{cc} 0 \ 1 \\ b \ d \\ \end{array}](2u_2 , 2v_2) \
\nonumber\\
&&\hskip-15mm 
+\varTheta[\begin{array}{cc} 1 \ 0 \\ b \ d \\ \end{array}](2u_1 , 2v_1) \ 
\varTheta[\begin{array}{cc} 1 \ 0 \\ b \ d \\ \end{array}](2u_2 , 2v_2) 
+\varTheta[\begin{array}{cc} 1 \ 1 \\ b \ d \\ \end{array}](2u_1 , 2v_1) \ 
\varTheta[\begin{array}{cc} 1 \ 1 \\ b \ d \\ \end{array}](2u_2 , 2v_2)  .
\label{2e7}
\end{eqnarray}
We give explicit formulae for each $(b,d)=(0,0), (0,1), (1,0), (1,1)$
of the form
\begin{eqnarray}
&&\hskip-15mm 1)\
\vartheta[\begin{array}{cc} 0 \ 0 \\ 0 \ 0 \\ \end{array}](u_1+u_2 ,v_1+v_2) \ 
\vartheta[\begin{array}{cc} 0 \ 0 \\ 0 \ 0 \\ \end{array}](u_1-u_2, v_1-v_2) 
\nonumber\\
&&\hskip-15mm 
=\varTheta[\begin{array}{cc} 0 \ 0 \\ 0 \ 0 \\ \end{array}](2u_1 , 2v_1) \  
\varTheta[\begin{array}{cc} 0  \ 0 \\ 0 \ 0 \\ \end{array}](2u_2 , 2v_2) 
+\varTheta[\begin{array}{cc} 0 \ 1 \\ 0 \ 0 \\ \end{array}](2u_1 , 2v_1) \ 
\varTheta[\begin{array}{cc} 0 \ 1 \\ 0 \ 0 \\ \end{array}](2u_2 , 2v_2) \
\nonumber\\
&&\hskip-15mm 
+\varTheta[\begin{array}{cc} 1 \ 0 \\ 0 \ 0 \\ \end{array}](2u_1 , 2v_1) \ 
\varTheta[\begin{array}{cc} 1 \ 0 \\ 0 \ 0 \\ \end{array}](2u_2 , 2v_2) 
+\varTheta[\begin{array}{cc} 1 \ 1 \\ 0 \ 0 \\ \end{array}](2u_1 , 2v_1) \ 
\varTheta[\begin{array}{cc} 1 \ 1 \\ 0 \ 0 \\ \end{array}](2u_2 , 2v_2)  ,
\label{2e8}
\end{eqnarray}
\begin{eqnarray}
&&\hskip-15mm 2)\
\vartheta[\begin{array}{cc} 0 \ 0 \\ 0 \ 1 \\ \end{array}](u_1+u_2 ,v_1+v_2) \ 
\vartheta[\begin{array}{cc} 0 \ 0 \\ 0 \ 0 \\ \end{array}](u_1-u_2, v_1-v_2) 
\nonumber\\
&&\hskip-15mm 
=\varTheta[\begin{array}{cc} 0 \ 0 \\ 0 \ 1 \\ \end{array}](2u_1 , 2v_1) \  
\varTheta[\begin{array}{cc} 0  \ 0 \\ 0 \ 1 \\ \end{array}](2u_2 , 2v_2) 
+\varTheta[\begin{array}{cc} 0 \ 1 \\ 0 \ 1 \\ \end{array}](2u_1 , 2v_1) \ 
\varTheta[\begin{array}{cc} 0 \ 1 \\0 \ 1 \\ \end{array}](2u_2 , 2v_2) \
\nonumber\\
&&\hskip-15mm 
+\varTheta[\begin{array}{cc} 1 \ 0 \\ 0 \ 1 \\ \end{array}](2u_1 , 2v_1) \ 
\varTheta[\begin{array}{cc} 1 \ 0 \\ 0 \ 1 \\ \end{array}](2u_2 , 2v_2) 
+\varTheta[\begin{array}{cc} 1 \ 1 \\ 0 \ 1 \\ \end{array}](2u_1 , 2v_1) \ 
\varTheta[\begin{array}{cc} 1 \ 1 \\ 0 \ 1 \\ \end{array}](2u_2 , 2v_2)  ,
\label{2e9}
\end{eqnarray}
\begin{eqnarray}
&&\hskip-15mm 3)\
\vartheta[\begin{array}{cc} 0 \ 0 \\ 1 \ 0 \\ \end{array}](u_1+u_2 ,v_1+v_2) \ 
\vartheta[\begin{array}{cc} 0 \ 0 \\ 0 \ 0 \\ \end{array}](u_1-u_2, v_1-v_2) 
\nonumber\\
&&\hskip-15mm 
=\varTheta[\begin{array}{cc} 0 \ 0 \\ 1 \ 0 \\ \end{array}](2u_1 , 2v_1) \  
\varTheta[\begin{array}{cc} 0  \ 0 \\ 1 \ 0 \\ \end{array}](2u_2 , 2v_2) 
+\varTheta[\begin{array}{cc} 0 \ 1 \\ 1 \ 0 \\ \end{array}](2u_1 , 2v_1) \ 
\varTheta[\begin{array}{cc} 0 \ 1 \\ 1 \ 0 \\ \end{array}](2u_2 , 2v_2)   
\nonumber\\
&&\hskip-15mm 
+\varTheta[\begin{array}{cc} 1 \ 0 \\ 1 \ 0 \\ \end{array}](2u_1 , 2v_1) \ 
\varTheta[\begin{array}{cc} 1 \ 0 \\ 1 \ 0 \\ \end{array}](2u_2 , 2v_2) 
+\varTheta[\begin{array}{cc} 1 \ 1 \\ 1 \ 0 \\ \end{array}](2u_1 , 2v_1) \ 
\varTheta[\begin{array}{cc} 1 \ 1 \\ 1 \ 0 \\ \end{array}](2u_2 , 2v_2)  ,
\label{2e10}
\end{eqnarray}
\begin{eqnarray}
&&\hskip-15mm 4)\
\vartheta[\begin{array}{cc} 0 \ 0 \\ 1 \ 1 \\ \end{array}](u_1+u_2 ,v_1+v_2) \ 
\vartheta[\begin{array}{cc} 0 \ 0 \\ 0 \ 0 \\ \end{array}](u_1-u_2, v_1-v_2) 
\nonumber\\
&&\hskip-15mm 
=\varTheta[\begin{array}{cc} 0 \ 0 \\ 1 \ 1 \\ \end{array}](2u_1 , 2v_1) \  
\varTheta[\begin{array}{cc} 0  \ 0 \\ 1 \ 1 \\ \end{array}](2u_2 , 2v_2) 
+\varTheta[\begin{array}{cc} 0 \ 1 \\ 1 \ 1 \\ \end{array}](2u_1 , 2v_1) \ 
\varTheta[\begin{array}{cc} 0 \ 1 \\ 1 \ 1 \\ \end{array}](2u_2 , 2v_2) \
\nonumber\\
&&\hskip-15mm 
+\varTheta[\begin{array}{cc} 1 \ 0 \\ 1 \ 1 \\ \end{array}](2u_1 , 2v_1) \ 
\varTheta[\begin{array}{cc} 1 \ 0 \\ 1 \ 1 \\ \end{array}](2u_2 , 2v_2) 
+\varTheta[\begin{array}{cc} 1 \ 1 \\ 1 \ 1 \\ \end{array}](2u_1 , 2v_1) \ 
\varTheta[\begin{array}{cc} 1 \ 1 \\ 1 \ 1 \\ \end{array}](2u_2 , 2v_2)  .
\label{2e11}
\end{eqnarray}
We eliminate  
$\displaystyle{\vartheta[\begin{array}{cc} 0 \ 0 \\ 0 \ 0 \\ \end{array}](u_1-u_2, v_1-v_2)}$
by taking the ratio between Eq.(\ref{2e8}) $\sim$ Eq.(\ref{2e11}), and we obtain three addition 
formulae. Detailed expressions are given in Appendix A by using Appendix B.\\
\\
\noindent
{\bf B) sector:\ $a=0$, $c=1$};\\
In this sector of Eq.(\ref{2e4}), addition formulae for
$\vartheta[\begin{array}{cc} 0 \ 1 \\ * \ * \\ \end{array}](u_1+u_2 ,v_1+v_2)$
are obtained by using following formulae
\begin{eqnarray}
&&\hskip -15mm
\vartheta[\begin{array}{cc} 0 \ 1 \\ b \ d \\ \end{array}](u_1+u_2 ,v_1+v_2) \ 
\vartheta[\begin{array}{cc} 0 \ 1 \\ 0 \ 0 \\ \end{array}](u_1-u_2, v_1-v_2) 
\nonumber\\
&&\hskip -15mm 
=\varTheta[\begin{array}{cc} 0 \ 1 \\ b \ d \\ \end{array}](2u_1 , 2v_1) \  
\varTheta[\begin{array}{cc} 0  \ 0 \\ b \ d \\ \end{array}](2u_2 , 2v_2) 
+\varTheta[\begin{array}{cc} 0 \ 0 \\ b \ d \\ \end{array}](2u_1 , 2v_1) \ 
\varTheta[\begin{array}{cc} 0 \ 1 \\ b \ d \\ \end{array}](2u_2 , 2v_2) \
\nonumber\\
&&\hskip -15mm
+\varTheta[\begin{array}{cc} 1 \ 1 \\ b \ d \\ \end{array}](2u_1 , 2v_1) \
\varTheta[\begin{array}{cc} 1 \ 0 \\ b \ d \\ \end{array}](2u_2 , 2v_2) 
+\varTheta[\begin{array}{cc} 1 \ 0 \\ b \ d \\ \end{array}](2u_1 , 2v_1) \ 
\varTheta[\begin{array}{cc} 1 \ 1 \\ b \ d \\ \end{array}](2u_2 , 2v_2)  .
\label{2e12}
\end{eqnarray}
We give explicit formulae for each $(b,d)=(0,0), (0,1), (1,0), (1,1)$
of the form
\begin{eqnarray}
&&\hskip -15mm 5)\
\vartheta[\begin{array}{cc} 0 \ 1 \\ 0 \ 0 \\ \end{array}](u_1+u_2 ,v_1+v_2) \ 
\vartheta[\begin{array}{cc} 0 \ 1 \\ 0 \ 0 \\ \end{array}](u_1-u_2, v_1-v_2) 
\nonumber\\
&&\hskip -15mm 
=\varTheta[\begin{array}{cc} 0 \ 1 \\ 0 \ 0 \\ \end{array}](2u_1 , 2v_1) \  
\varTheta[\begin{array}{cc} 0  \ 0 \\ 0 \ 0 \\ \end{array}](2u_2 , 2v_2) 
+\varTheta[\begin{array}{cc} 0 \ 0 \\ 0 \ 0 \\ \end{array}](2u_1 , 2v_1) \ 
\varTheta[\begin{array}{cc} 0 \ 1 \\ 0 \ 0 \\ \end{array}](2u_2 , 2v_2) \
\nonumber\\
&&\hskip -15mm
+\varTheta[\begin{array}{cc} 1 \ 1 \\ 0 \ 0 \\ \end{array}](2u_1 , 2v_1) \
\varTheta[\begin{array}{cc} 1 \ 0 \\ 0 \ 0 \\ \end{array}](2u_2 , 2v_2) 
+\varTheta[\begin{array}{cc} 1 \ 0 \\ 0 \ 0 \\ \end{array}](2u_1 , 2v_1) \ 
\varTheta[\begin{array}{cc} 1 \ 1 \\ 0 \ 0 \\ \end{array}](2u_2 , 2v_2)  ,
\label{2e13}
\end{eqnarray}
\begin{eqnarray}
&&\hskip -15mm  6)\
\vartheta[\begin{array}{cc} 0 \ 1 \\ 0 \ 1 \\ \end{array}](u_1+u_2 ,v_1+v_2) \ 
\vartheta[\begin{array}{cc} 0 \ 1 \\ 0 \ 0 \\ \end{array}](u_1-u_2, v_1-v_2) 
\nonumber\\
&&\hskip -15mm 
=\varTheta[\begin{array}{cc} 0 \ 1 \\ 0 \ 1 \\ \end{array}](2u_1 , 2v_1) \  
\varTheta[\begin{array}{cc} 0  \ 0 \\ 0 \ 1 \\ \end{array}](2u_2 , 2v_2) 
+\varTheta[\begin{array}{cc} 0 \ 0 \\ 0 \ 1 \\ \end{array}](2u_1 , 2v_1) \ 
\varTheta[\begin{array}{cc} 0 \ 1 \\ 0 \ 1 \\ \end{array}](2u_2 , 2v_2) \
\nonumber\\
&&\hskip -15mm
+\varTheta[\begin{array}{cc} 1 \ 1 \\ 0 \ 1 \\ \end{array}](2u_1 , 2v_1) \
\varTheta[\begin{array}{cc} 1 \ 0 \\ 0 \ 1 \\ \end{array}](2u_2 , 2v_2) 
+\varTheta[\begin{array}{cc} 1 \ 0 \\ 0 \ 1 \\ \end{array}](2u_1 , 2v_1) \ 
\varTheta[\begin{array}{cc} 1 \ 1 \\ 0 \ 1 \\ \end{array}](2u_2 , 2v_2)  ,
\label{2e14}
\end{eqnarray}
\begin{eqnarray}
&&\hskip -15mm 7)\
\vartheta[\begin{array}{cc} 0 \ 1 \\ 1 \ 0 \\ \end{array}](u_1+u_2 ,v_1+v_2) \ 
\vartheta[\begin{array}{cc} 0 \ 1 \\ 0 \ 0 \\ \end{array}](u_1-u_2, v_1-v_2) 
\nonumber\\
&&\hskip -15mm 
=\varTheta[\begin{array}{cc} 0 \ 1 \\ 1 \ 0 \\ \end{array}](2u_1 , 2v_1) \  
\varTheta[\begin{array}{cc} 0  \ 0 \\ 1 \ 0 \\ \end{array}](2u_2 , 2v_2) 
+\varTheta[\begin{array}{cc} 0 \ 0 \\ 1 \ 0 \\ \end{array}](2u_1 , 2v_1) \ 
\varTheta[\begin{array}{cc} 0 \ 1 \\ 1 \ 0 \\ \end{array}](2u_2 , 2v_2) \
\nonumber\\
&&\hskip -15mm
+\varTheta[\begin{array}{cc} 1 \ 1 \\ 1 \ 0 \\ \end{array}](2u_1 , 2v_1) \
\varTheta[\begin{array}{cc} 1 \ 0 \\ 1 \ 0 \\ \end{array}](2u_2 , 2v_2) 
+\varTheta[\begin{array}{cc} 1 \ 0 \\ 1 \ 0 \\ \end{array}](2u_1 , 2v_1) \ 
\varTheta[\begin{array}{cc} 1 \ 1 \\ 1 \ 0 \\ \end{array}](2u_2 , 2v_2)  ,
\label{2e15}
\end{eqnarray}
\begin{eqnarray}
&&\hskip -15mm 8)\
\vartheta[\begin{array}{cc} 0 \ 1 \\ 1 \ 1 \\ \end{array}](u_1+u_2 ,v_1+v_2) \ 
\vartheta[\begin{array}{cc} 0 \ 1 \\ 0 \ 0 \\ \end{array}](u_1-u_2, v_1-v_2) 
\nonumber\\
&&\hskip -15mm 
=\varTheta[\begin{array}{cc} 0 \ 1 \\ 1 \ 1 \\ \end{array}](2u_1 , 2v_1) \  
\varTheta[\begin{array}{cc} 0  \ 0 \\ 1 \ 1 \\ \end{array}](2u_2 , 2v_2) 
+\varTheta[\begin{array}{cc} 0 \ 0 \\ 1 \ 1 \\ \end{array}](2u_1 , 2v_1) \ 
\varTheta[\begin{array}{cc} 0 \ 1 \\ 1 \ 1 \\ \end{array}](2u_2 , 2v_2) \
\nonumber\\
&&\hskip -15mm
+\varTheta[\begin{array}{cc} 1 \ 1 \\ 1 \ 1 \\ \end{array}](2u_1 , 2v_1) \
\varTheta[\begin{array}{cc} 1 \ 0 \\ 1 \ 1 \\ \end{array}](2u_2 , 2v_2) 
+\varTheta[\begin{array}{cc} 1 \ 0 \\ 1 \ 1 \\ \end{array}](2u_1 , 2v_1) \ 
\varTheta[\begin{array}{cc} 1 \ 1 \\ 1 \ 1 \\ \end{array}](2u_2 , 2v_2)  .
\label{2e16}
\end{eqnarray}
We eliminate  $\displaystyle{\vartheta[\begin{array}{cc} 0 \ 1 \\ 0 \ 0 \\ \end{array}](u_1-u_2, v_1-v_2)}$
by taking the ratio between Eq.(\ref{2e13}) $\sim$ Eq.(\ref{2e16}), and we obtain three addition 
formulae. Detailed expressions are given in Appendix A by using Appendix B.\\
%
\\
\noindent
{\bf C) sector:\ $a=1$, $c=0$}; \\
In this sector of Eq.(\ref{2e4}), addition formulae for
$\vartheta[\begin{array}{cc} 1 \ 0 \\ * \ * \\ \end{array}](u_1+u_2 ,v_1+v_2)$
are obtained by using following formulae
\begin{eqnarray}
&&\hskip-15mm 
\vartheta[\begin{array}{cc} 1 \ 0 \\ b \ d \\ \end{array}](u_1+u_2 ,v_1+v_2) \ 
\vartheta[\begin{array}{cc} 1 \ 0 \\ 0 \ 0 \\ \end{array}](u_1-u_2, v_1-v_2) 
\nonumber\\
&&\hskip-15mm 
=\varTheta[\begin{array}{cc} 1 \ 0 \\ b \ d \\ \end{array}](2u_1 , 2v_1) \  
\varTheta[\begin{array}{cc} 0  \ 0 \\ b \ d \\ \end{array}](2u_2 , 2v_2) 
+\varTheta[\begin{array}{cc} 1 \ 1 \\ b \ d \\ \end{array}](2u_1 , 2v_1) \ 
\varTheta[\begin{array}{cc} 0 \ 1 \\ b \ d \\ \end{array}](2u_2 , 2v_2) \
\nonumber\\
&&\hskip-15mm 
+\varTheta[\begin{array}{cc} 0 \ 0 \\ b \ d \\ \end{array}](2u_1 , 2v_1) \
\varTheta[\begin{array}{cc} 1 \ 0 \\ b \ d \\ \end{array}](2u_2 , 2v_2) 
+\varTheta[\begin{array}{cc} 0 \ 1 \\ b \ d \\ \end{array}](2u_1 , 2v_1) \ 
\varTheta[\begin{array}{cc} 1 \ 1 \\ b \ d \\ \end{array}](2u_2 , 2v_2)  .
\label{2e17}
\end{eqnarray}
We give explicit formulae for each $(b,d)=(0,0), (0,1), (1,0), (1,1)$
of the form
\begin{eqnarray}
&&\hskip-15mm 9)\
\vartheta[\begin{array}{cc} 1 \ 0 \\ 0 \ 0 \\ \end{array}](u_1+u_2 ,v_1+v_2) \ 
\vartheta[\begin{array}{cc} 1 \ 0 \\ 0 \ 0 \\ \end{array}](u_1-u_2, v_1-v_2) 
\nonumber\\
&&\hskip-15mm 
=\varTheta[\begin{array}{cc} 1 \ 0 \\ 0 \ 0 \\ \end{array}](2u_1 , 2v_1) \  
\varTheta[\begin{array}{cc} 0  \ 0 \\ 0 \ 0 \\ \end{array}](2u_2 , 2v_2) 
+\varTheta[\begin{array}{cc} 1 \ 1 \\ 0 \ 0 \\ \end{array}](2u_1 , 2v_1) \ 
\varTheta[\begin{array}{cc} 0 \ 1 \\ 0 \ 0 \\ \end{array}](2u_2 , 2v_2) \
\nonumber\\
&&\hskip-15mm 
+\varTheta[\begin{array}{cc} 0 \ 0 \\ 0 \ 0 \\ \end{array}](2u_1 , 2v_1) \
\varTheta[\begin{array}{cc} 1 \ 0 \\ 0 \ 0 \\ \end{array}](2u_2 , 2v_2) 
+\varTheta[\begin{array}{cc} 0 \ 1 \\ 0 \ 0 \\ \end{array}](2u_1 , 2v_1) \ 
\varTheta[\begin{array}{cc} 1 \ 1 \\ 0 \ 0 \\ \end{array}](2u_2 , 2v_2) ,
\label{2e18}
\end{eqnarray}
\begin{eqnarray}
&&\hskip-15mm 10)\ 
\vartheta[\begin{array}{cc} 1 \ 0 \\ 0 \ 1 \\ \end{array}](u_1+u_2 ,v_1+v_2) \ 
\vartheta[\begin{array}{cc} 1 \ 0 \\ 0 \ 0 \\ \end{array}](u_1-u_2, v_1-v_2) 
\nonumber\\
&&\hskip-15mm 
=\varTheta[\begin{array}{cc} 1 \ 0 \\ 0 \ 1 \\ \end{array}](2u_1 , 2v_1) \  
\varTheta[\begin{array}{cc} 0  \ 0 \\ 0 \ 1 \\ \end{array}](2u_2 , 2v_2) 
+\varTheta[\begin{array}{cc} 1 \ 1 \\ 0 \ 1 \\ \end{array}](2u_1 , 2v_1) \ 
\varTheta[\begin{array}{cc} 0 \ 1 \\ 0 \ 1 \\ \end{array}](2u_2 , 2v_2) \
\nonumber\\
&&\hskip-15mm 
+\varTheta[\begin{array}{cc} 0 \ 0 \\ 0 \ 1 \\ \end{array}](2u_1 , 2v_1) \
\varTheta[\begin{array}{cc} 1 \ 0 \\ 0 \ 1 \\ \end{array}](2u_2 , 2v_2) 
+\varTheta[\begin{array}{cc} 0 \ 1 \\ 0 \ 1 \\ \end{array}](2u_1 , 2v_1) \ 
\varTheta[\begin{array}{cc} 1 \ 1 \\ 0 \ 1 \\ \end{array}](2u_2 , 2v_2) ,
\label{2e19}
\end{eqnarray}
\begin{eqnarray}
&&\hskip-15mm 11)\
\vartheta[\begin{array}{cc} 1 \ 0 \\ 1 \ 0 \\ \end{array}](u_1+u_2 ,v_1+v_2) \ 
\vartheta[\begin{array}{cc} 1 \ 0 \\ 0 \ 0 \\ \end{array}](u_1-u_2, v_1-v_2) 
\nonumber\\
&&\hskip-15mm 
=\varTheta[\begin{array}{cc} 1 \ 0 \\ 1 \ 0 \\ \end{array}](2u_1 , 2v_1) \  
\varTheta[\begin{array}{cc} 0  \ 0 \\ 1 \ 0 \\ \end{array}](2u_2 , 2v_2) 
+\varTheta[\begin{array}{cc} 1 \ 1 \\ 1 \ 0 \\ \end{array}](2u_1 , 2v_1) \ 
\varTheta[\begin{array}{cc} 0 \ 1 \\ 1 \ 0 \\ \end{array}](2u_2 , 2v_2) \
\nonumber\\
&&\hskip-15mm 
+\varTheta[\begin{array}{cc} 0 \ 0 \\ 1 \ 0 \\ \end{array}](2u_1 , 2v_1) \
\varTheta[\begin{array}{cc} 1 \ 0 \\ 1 \ 0 \\ \end{array}](2u_2 , 2v_2) 
+\varTheta[\begin{array}{cc} 0 \ 1 \\ 1 \ 0 \\ \end{array}](2u_1 , 2v_1) \ 
\varTheta[\begin{array}{cc} 1 \ 1 \\ 1 \ 0 \\ \end{array}](2u_2 , 2v_2) ,
\label{2e20}
\end{eqnarray}
\begin{eqnarray}
&&\hskip-15mm 12)\
\vartheta[\begin{array}{cc} 1 \ 0 \\ 1 \ 1 \\ \end{array}](u_1+u_2 ,v_1+v_2) \ 
\vartheta[\begin{array}{cc} 1 \ 0 \\ 0 \ 0 \\ \end{array}](u_1-u_2, v_1-v_2) 
\nonumber\\
&&\hskip-15mm 
=\varTheta[\begin{array}{cc} 1 \ 0 \\ 1 \ 1 \\ \end{array}](2u_1 , 2v_1) \  
\varTheta[\begin{array}{cc} 0  \ 0 \\ 1 \ 1 \\ \end{array}](2u_2 , 2v_2) 
+\varTheta[\begin{array}{cc} 1 \ 1 \\ 1 \ 1 \\ \end{array}](2u_1 , 2v_1) \ 
\varTheta[\begin{array}{cc} 0 \ 1 \\ 1 \ 1 \\ \end{array}](2u_2 , 2v_2) \
\nonumber\\
&&\hskip-15mm 
+\varTheta[\begin{array}{cc} 0 \ 0 \\ 1 \ 1 \\ \end{array}](2u_1 , 2v_1) \
\varTheta[\begin{array}{cc} 1 \ 0 \\ 1 \ 1 \\ \end{array}](2u_2 , 2v_2) 
+\varTheta[\begin{array}{cc} 0 \ 1 \\ 1 \ 1 \\ \end{array}](2u_1 , 2v_1) \ 
\varTheta[\begin{array}{cc} 1 \ 1 \\ 1 \ 1 \\ \end{array}](2u_2 , 2v_2) .
\label{2e21}
\end{eqnarray}
We eliminate  $\displaystyle{\vartheta[\begin{array}{cc} 1 \ 0 \\ 0 \ 0 \\ \end{array}](u_1-u_2, v_1-v_2)}$
by taking the ratio between Eq.(\ref{2e18}) $\sim$ Eq.(\ref{2e21}), and 
we obtain three addition 
formulae. Detailed expressions are given in Appendix A by using Appendix B.\\
%
\\
\noindent
{\bf D) sector\ $a=1$, $c=1$};\\
In this sector of Eq.(\ref{2e4}), addition formulae for
$\vartheta[\begin{array}{cc} 1 \ 1 \\ * \ * \\ \end{array}](u_1+u_2 ,v_1+v_2)$
are obtained by using following formulae
\begin{eqnarray}
&&\hskip-15mm 
\vartheta[\begin{array}{cc} 1 \ 1 \\ b \ d \\ \end{array}](u_1+u_2 ,v_1+v_2) \ 
\vartheta[\begin{array}{cc} 1 \ 1 \\ 0 \ 0 \\ \end{array}](u_1-u_2, v_1-v_2) 
\nonumber\\
&&\hskip-15mm 
=\varTheta[\begin{array}{cc} 1 \ 1 \\ b \ d \\ \end{array}](2u_1 , 2v_1) \  
\varTheta[\begin{array}{cc} 0  \ 0 \\ b \ d \\ \end{array}](2u_2 , 2v_2) 
+\varTheta[\begin{array}{cc} 1 \ 0 \\ b \ d \\ \end{array}](2u_1 , 2v_1) \ 
\varTheta[\begin{array}{cc} 0 \ 1 \\ b \ d \\ \end{array}](2u_2 , 2v_2) \
\nonumber\\
&&\hskip-15mm 
+\varTheta[\begin{array}{cc} 0 \ 1 \\ b \ d \\ \end{array}](2u_1 , 2v_1) \ 
\varTheta[\begin{array}{cc} 1 \ 0 \\ b \ d \\ \end{array}](2u_2 , 2v_2) 
+\varTheta[\begin{array}{cc} 0 \ 0 \\ b \ d \\ \end{array}](2u_1 , 2v_1) \ 
\varTheta[\begin{array}{cc} 1 \ 1 \\ b \ d \\ \end{array}](2u_2 , 2v_2) .
\label{2e22}
\end{eqnarray}
We give explicit formulae for each $(b,d)=(0,0), (0,1), (1,0), (1,1)$
of the form
\begin{eqnarray}
&&\hskip-15mm 13)\
\vartheta[\begin{array}{cc} 1 \ 1 \\ 0 \ 0 \\ \end{array}](u_1+u_2 ,v_1+v_2) \ 
\vartheta[\begin{array}{cc} 1 \ 1 \\ 0 \ 0 \\ \end{array}](u_1-u_2, v_1-v_2) 
\nonumber\\
&&\hskip-15mm 
=\varTheta[\begin{array}{cc} 1 \ 1 \\ 0 \ 0 \\ \end{array}](2u_1 , 2v_1) \  
\varTheta[\begin{array}{cc} 0  \ 0 \\ 0 \ 0 \\ \end{array}](2u_2 , 2v_2) 
+\varTheta[\begin{array}{cc} 1 \ 0 \\ 0 \ 0 \\ \end{array}](2u_1 , 2v_1) \ 
\varTheta[\begin{array}{cc} 0 \ 1 \\ 0 \ 0 \\ \end{array}](2u_2 , 2v_2) \
\nonumber\\
&&\hskip-15mm 
+\varTheta[\begin{array}{cc} 0 \ 1 \\ 0 \ 0 \\ \end{array}](2u_1 , 2v_1) \ 
\varTheta[\begin{array}{cc} 1 \ 0 \\ 0 \ 0 \\ \end{array}](2u_2 , 2v_2) 
+\varTheta[\begin{array}{cc} 0 \ 0 \\ 0 \ 0 \\ \end{array}](2u_1 , 2v_1) \ 
\varTheta[\begin{array}{cc} 1 \ 1 \\ 0 \ 0 \\ \end{array}](2u_2 , 2v_2) ,
\label{2e23}
\end{eqnarray}
\begin{eqnarray}
&&\hskip-15mm 14)\
\vartheta[\begin{array}{cc} 1 \ 1 \\ 0 \ 1 \\ \end{array}](u_1+u_2 ,v_1+v_2) \ 
\vartheta[\begin{array}{cc} 1 \ 1 \\ 0  0 \\ \end{array}](u_1-u_2, v_1-v_2) 
\nonumber\\
&&\hskip-15mm 
=\varTheta[\begin{array}{cc} 1 \ 1 \\ 0 \ 1 \\ \end{array}](2u_1 , 2v_1) \  
\varTheta[\begin{array}{cc} 0  \ 0 \\ 0 \ 1 \\ \end{array}](2u_2 , 2v_2) 
+\varTheta[\begin{array}{cc} 1 \ 0 \\ 0 \ 1 \\ \end{array}](2u_1 , 2v_1) \ 
\varTheta[\begin{array}{cc} 0 \ 1 \\ 0 \ 1 \\ \end{array}](2u_2 , 2v_2) \
\nonumber\\
&&\hskip-15mm 
+\varTheta[\begin{array}{cc} 0 \ 1 \\ 0 \ 1 \\ \end{array}](2u_1 , 2v_1) \ 
\varTheta[\begin{array}{cc} 1 \ 0 \\ 0 \ 1 \\ \end{array}](2u_2 , 2v_2) 
+\varTheta[\begin{array}{cc} 0 \ 0 \\ 0 \ 1 \\ \end{array}](2u_1 , 2v_1) \ 
\varTheta[\begin{array}{cc} 1 \ 1 \\ 0 \ 1 \\ \end{array}](2u_2 , 2v_2) ,
\label{2e24}
\end{eqnarray}
\begin{eqnarray}
&&\hskip-15mm 15)\
\vartheta[\begin{array}{cc} 1 \ 1 \\ 1 \ 0 \\ \end{array}](u_1+u_2 ,v_1+v_2) \ 
\vartheta[\begin{array}{cc} 1 \ 1 \\ 0 \ 0 \\ \end{array}](u_1-u_2, v_1-v_2) 
\nonumber\\
&&\hskip-15mm 
=\varTheta[\begin{array}{cc} 1 \ 1 \\ 1 \ 0 \\ \end{array}](2u_1 , 2v_1) \  
\varTheta[\begin{array}{cc} 0  \ 0 \\ 1 \ 0 \\ \end{array}](2u_2 , 2v_2) 
+\varTheta[\begin{array}{cc} 1 \ 0 \\ 1 \ 0 \\ \end{array}](2u_1 , 2v_1) \ 
\varTheta[\begin{array}{cc} 0 \ 1 \\ 1 \ 0 \\ \end{array}](2u_2 , 2v_2) \
\nonumber\\
&&\hskip-15mm 
+\varTheta[\begin{array}{cc} 0 \ 1 \\ 1 \ 0 \\ \end{array}](2u_1 , 2v_1) \ 
\varTheta[\begin{array}{cc} 1 \ 0 \\ 1 \ 0 \\ \end{array}](2u_2 , 2v_2) 
+\varTheta[\begin{array}{cc} 0 \ 0 \\ 1 \ 0 \\ \end{array}](2u_1 , 2v_1) \ 
\varTheta[\begin{array}{cc} 1 \ 1 \\ 1 \ 0 \\ \end{array}](2u_2 , 2v_2) ,
\label{2e25}
\end{eqnarray}
\begin{eqnarray}
&&\hskip-15mm 16)\
\vartheta[\begin{array}{cc} 1 \ 1 \\ 1 \ 1 \\ \end{array}](u_1+u_2 ,v_1+v_2) \ 
\vartheta[\begin{array}{cc} 1 \ 1 \\ 0 \ 0 \\ \end{array}](u_1-u_2, v_1-v_2) 
\nonumber\\
&&\hskip-15mm 
=\varTheta[\begin{array}{cc} 1 \ 1 \\ 1 \ 1 \\ \end{array}](2u_1 , 2v_1) \  
\varTheta[\begin{array}{cc} 0  \ 0 \\ 1 \ 1 \\ \end{array}](2u_2 , 2v_2) 
+\varTheta[\begin{array}{cc} 1 \ 0 \\ 1 \ 1 \\ \end{array}](2u_1 , 2v_1) \ 
\varTheta[\begin{array}{cc} 0 \ 1 \\ 1 \ 1 \\ \end{array}](2u_2 , 2v_2) \
\nonumber\\
&&\hskip-15mm 
+\varTheta[\begin{array}{cc} 0 \ 1 \\ 1 \ 1 \\ \end{array}](2u_1 , 2v_1) \ 
\varTheta[\begin{array}{cc} 1 \ 0 \\ 1 \ 1 \\ \end{array}](2u_2 , 2v_2) 
+\varTheta[\begin{array}{cc} 0 \ 0 \\ 1 \ 1 \\ \end{array}](2u_1 , 2v_1) \ 
\varTheta[\begin{array}{cc} 1 \ 1 \\ 1 \ 1 \\ \end{array}](2u_2 , 2v_2)  .
\label{2e26}
\end{eqnarray}
We eliminate  $\displaystyle{\vartheta[\begin{array}{cc} 1 \ 1 \\ 0 \ 0 \\ \end{array}](u_1-u_2, v_1-v_2)}$
by taking the ratio between Eq.(\ref{2e23}) $\sim$ Eq.(\ref{2e26}), and we obtain three addition 
formulae. Detailed expressions are given in Appendix A by using Appendix B.\\
%
\subsection{Functional Relations}
In order to give addition formulae, we must express 
$\varTheta\left[\begin{array}{cc} * \ * \\ * \ * \\ \end{array}\right](2u, 2v)$
by $\vartheta\left[\begin{array}{cc} * \ * \\ * \ * \\ \end{array}\right](u,v)$. \\
\\
{\bf Functional relations 1)}\\ 
We obtain
$\varTheta\left[\begin{array}{cc} * \ * \\ 0 \ 0 \\ \end{array}\right](2u, 2v)$ 
 by putting $a=0$, $c=0$ in functional relations(I) of Eq.(\ref{2e5}) of the form
\begin{eqnarray}
&&\hskip-15mm \vartheta\left[\begin{array}{cc} 0 \ 0 \\ b \ d \\ \end{array}\right]^2(u ,v) \ 
\nonumber\\
&&\hskip-15mm 
=\varTheta\left[\begin{array}{cc} 0 \ 0 \\ 2b \ 2d \\ \end{array}\right] (2u , 2v) \ \ 
\varTheta\left[\begin{array}{cc} 0 \ 0 \\ 0 \ 0 \\ \end{array}\right] (0 ,0) 
+\varTheta\left[\begin{array}{cc} 0 \ 1 \\ 2b \ 2d \\ \end{array}\right] (2u , 2v) \ \
\varTheta\left[\begin{array}{cc} 0 \ 1 \\ 0 \ 0 \\ \end{array}\right](0, 0) \
\nonumber\\
&&\hskip-15mm +\varTheta\left[\begin{array}{cc} 1 \ 0 \\ 2b \ 2d \\ \end{array}\right] (2u , 2v) \ \
\varTheta\left[\begin{array}{cc} 1 \ 0 \\ 0 \ 0 \\ \end{array}\right](0 , 0) 
+\varTheta\left[\begin{array}{cc} 1 \ 1 \\ 2b \ 2d \\ \end{array}\right](2u , 2v) \ \
\varTheta\left[\begin{array}{cc} 1 \ 1 \\ 0 \ 0 \\ \end{array}\right](0 , 0) ,
\label{2e27}
\end{eqnarray}
We give the explicit formula for each $(b,d)=(0,0), (0,1), (1,0), (1,1)$
in the following matrix form
\begin{eqnarray}
  &&\hspace{-20mm} 
  \left( %
 \right)   ,
 \label{2e28}
\end{eqnarray}
Here the Riemann matrix, which is the self adjoint orthogonal matrix, naturally 
appears. 
We inversely solve the above in the following matrix form
\begin{eqnarray}
  &&\hspace{-20mm} 
  \left( %
] (0 , 0) $
are given by  $\vartheta[\begin{array}{cc} 0 \ 0 \\ * \ * \\ \end{array}] (0 , 0) $
in the next \S 2.4.   \\
\\
{\bf Functional relations 2)} \\
We obtain
$\varTheta\left[\begin{array}{cc} * \ * \\ 0 \ 1 \\ \end{array}\right](2u, 2v)$ 
by putting $b=0$, $d=1$ in function relations(II) of Eq.(\ref{2e6}) of the form
\begin{eqnarray}
&&\hskip-15mm \vartheta[\begin{array}{cc} a \ c \\ 0 \ 1 \\ \end{array}](u , v) \ 
\vartheta[\begin{array}{cc} a \ c \\ 0 \ 0 \\ \end{array}](u, v) 
\nonumber\\
&&\hskip-15mm
=\varTheta[\begin{array}{cc} a \ c \\ 0 \ 1 \\ \end{array}](2u , 2v) \  
\varTheta[\begin{array}{cc} 0 \ 0 \\ 0 \ 1 \\ \end{array}](0 , 0) 
+\varTheta[\begin{array}{cc} a+1 &\hspace{-2mm}  c \\ 0 &\hspace{-2mm}  1 \\ \end{array}](2u , 2v) \ 
\varTheta[\begin{array}{cc} 1 \ 0 \\ 0 \ 1 \\ \end{array}](0 , 0) ,
\label{2e30}
\end{eqnarray}
where we use $\displaystyle{\varTheta[\begin{array}{cc} 0 \ 1 \\ 0 \ 1 \\ \end{array}](0 , 0) =0}$
and $\displaystyle{\varTheta[\begin{array}{cc} 1 \ 1 \\ 0 \ 1 \\ \end{array}](0 , 0) =0}$.
We give explicit formulae for each $(a,c)=(0,0), (0,1), (1,0), (1,1)$
in the following matrix form
\begin{eqnarray}
  &&\hspace{-20mm} 
  \left( %
](0 ,0)}$.
Then we inversely solve the above in the following matrix form
\begin{eqnarray}
  &&\hspace{-30mm} 
  \left( %
\right](2u, 2v)$ 
 by putting $b=1$, $d=0$ in functional relations(II).
Using the left-right symmetry, we have the following matrix form from Eq.(\ref{2e32}) 
\begin{eqnarray}
  &&\hspace{-30mm} 
  \left( %
\right](2u, 2v)$ 
 by putting $b=1$, $d=1$ in the function relation(II) of the matrix form
\begin{eqnarray}
  &&\hspace{-20mm} 
  \left( %
](0 ,0)}$.
Then we inversely solve the above in the following matrix form
\begin{eqnarray}
  &&\hspace{-30mm} 
  \left( %
 \right)   .
 \label{2e35}
\end{eqnarray}
In this way, from Eq.(\ref{2e29}), Eq.(\ref{2e32}), Eq.(\ref{2e33}), Eq.(\ref{2e35}), we 
have expressions 
$\varTheta\left[\begin{array}{cc} * \ * \\ * \ * \\ \end{array}\right] (2 u , 2 v)$ by
$\vartheta\left[\begin{array}{cc} * \ * \\ * \ * \\ \end{array}\right] (u , v)$.
Explicit expressions are given in Appendix C.\\
\subsection{Constants Relations}
Finally, in order to give addition formulae, we must express 
$\varTheta\left[\begin{array}{cc} * \ * \\ * \ * \\ \end{array}\right](0, 0)$
by $\vartheta\left[\begin{array}{cc} * \ * \\ * \ * \\ \end{array}\right](0, 0)$. \\
\\
{\bf Constants relations 1)}\\
We obtain  $\varTheta[\begin{array}{cc} * \ * \\ 0 \ 0 \\ \end{array}]^2(0 ,0)$
by putting $u=0$, $v=0$ in Eq.(\ref{2e29}) in the form
\begin{eqnarray}
  &&\hspace{-20mm} 
  \left( %
]^2(0 ,0)$
by putting $u=0$, $v=0$ in Eq.(\ref{2e31}). From the first and the third component, we have
\begin{eqnarray}
&&\hskip-15mm  \vartheta\left[\begin{array}{cc} 0 \ 0 \\ 0 \ 1 \\ \end{array}\right] (0 , 0) \ 
 \vartheta\left[\begin{array}{cc} 0 \ 0 \\ 0 \ 0 \\ \end{array}\right] (0 , 0)  
=\alpha \varTheta\left[\begin{array}{cc} 0 \ 0 \\ 0 \ 1 \\ \end{array}\right]  (0 , 0) \  
+\beta  \varTheta\left[\begin{array}{cc} 1 \ 0 \\ 0 \ 1 \\ \end{array}\right] (0 , 0)
=\alpha^2+\beta^2, \nonumber\\
&&\hskip-15mm  \vartheta\left[\begin{array}{cc} 1 \ 0 \\ 0 \ 1 \\ \end{array}\right] (0 , 0) \ 
 \vartheta\left[\begin{array}{cc} 1 \ 0 \\ 0 \ 0 \\ \end{array}\right] (0 , 0)  
=\beta \varTheta\left[\begin{array}{cc} 0 \ 0 \\ 0 \ 1 \\ \end{array}\right]  (0 , 0) \  
+\alpha  \varTheta\left[\begin{array}{cc} 1 \ 0 \\ 0 \ 1 \\ \end{array}\right] (0 , 0)
=2 \alpha \beta   . 
\label{2e37}
\end{eqnarray}
This gives
\begin{eqnarray}
&&\hskip-15mm  (\alpha\pm \beta)^2=
\Big(\varTheta\left[\begin{array}{cc} 0 \ 0 \\ 0 \ 1 \\ \end{array}\right] (0 , 0) \ 
 \pm \varTheta\left[\begin{array}{cc} 1 \ 0 \\ 0 \ 1 \\ \end{array}\right] (0 , 0) \Big)^2 \nonumber\\
&&\hskip-15mm 
=\vartheta\left[\begin{array}{cc} 0 \ 0 \\ 0 \ 1 \\ \end{array}\right]  (0 , 0) \  
 \vartheta\left[\begin{array}{cc} 0 \ 0 \\ 0 \ 0 \\ \end{array}\right] (0 , 0) \ \
\pm\vartheta\left[\begin{array}{cc} 1 \ 0 \\ 0 \ 1 \\ \end{array}\right] (0 , 0) \ \
\vartheta\left[\begin{array}{cc} 1 \ 0 \\ 0 \ 0 \\ \end{array}\right] (0 , 0)    .
\label{2e38}
\end{eqnarray}
The others are  $\varTheta\left[\begin{array}{cc} 0 \ 1 \\ 0 \ 1 \\ \end{array}\right] (0 , 0)=0$,
 $\varTheta\left[\begin{array}{cc} 1 \ 1 \\ 0 \ 1 \\ \end{array}\right] (0 , 0)=0$ because these 
are odd functions. \\
%
\\
{\bf Constants relations 3)}\\
We obtain  $\varTheta[\begin{array}{cc} * \ * \\ 1 \ 0 \\ \end{array}]^2(0 ,0)$
by putting $u=0$, $v=0$ in Eq.(\ref{2e33}). From the first and the third component, 
which is obtained by exchanging the left-right index of Eq.(\ref{2e38}) of the form
\begin{eqnarray}
&&\hskip-15mm  \Big(\varTheta\left[\begin{array}{cc} 0 \ 0 \\ 1 \ 0 \\ \end{array}\right] (0 , 0) \ 
 \pm \varTheta\left[\begin{array}{cc} 0 \ 1 \\ 1 \ 0 \\ \end{array}\right] (0 , 0) \Big)^2 \nonumber\\
&&\hskip-15mm 
=\vartheta\left[\begin{array}{cc} 0 \ 0 \\ 1 \ 0 \\ \end{array}\right]  (0 , 0) \  
 \vartheta\left[\begin{array}{cc} 0 \ 0 \\ 0 \ 0 \\ \end{array}\right] (0 , 0) \ \
\pm\vartheta\left[\begin{array}{cc} 0 \ 1 \\ 1 \ 0 \\ \end{array}\right] (0 , 0) \ \
\vartheta\left[\begin{array}{cc} 0 \ 1 \\ 0 \ 0 \\ \end{array}\right] (0 , 0)  .
\label{2e39}
\end{eqnarray} 
\\
{\bf Constants relations 4)}\\
We obtain  $\varTheta[\begin{array}{cc} * \ * \\ 1 \ 1 \\ \end{array}]^2(0 ,0)$
by putting $u=0$, $v=0$ in Eq.(\ref{2e35}). From the first and the fourth components, we have
\begin{eqnarray}
&&\hskip-15mm  \vartheta\left[\begin{array}{cc} 0 \ 0 \\ 1 \ 1 \\ \end{array}\right] (0 , 0) \ 
 \vartheta\left[\begin{array}{cc} 0 \ 0 \\ 0 \ 0 \\ \end{array}\right] (0 , 0)  
=\xi \varTheta\left[\begin{array}{cc} 0 \ 0 \\ 1 \ 1 \\ \end{array}\right]  (0 , 0) \  
+\zeta  \varTheta\left[\begin{array}{cc} 1 \ 1 \\ 1 \ 1 \\ \end{array}\right] (0 , 0)
=\xi^2+\zeta^2, \nonumber\\
&&\hskip-15mm  \vartheta\left[\begin{array}{cc} 1 \ 1 \\ 1 \ 1 \\ \end{array}\right] (0 , 0) \ 
 \vartheta\left[\begin{array}{cc} 1 \ 1 \\ 0 \ 0 \\ \end{array}\right] (0 , 0)  
=\zeta \varTheta\left[\begin{array}{cc} 0 \ 0 \\ 1 \ 1 \\ \end{array}\right]  (0 , 0) \  
+\xi  \varTheta\left[\begin{array}{cc} 1 \ 1 \\ 1 \ 1 \\ \end{array}\right] (0 , 0)
=2 \xi \zeta  . 
\label{2e40}
\end{eqnarray}
This gives 
\begin{eqnarray}
&&\hskip-15mm  (\xi\pm \zeta)^2=
\Big(\varTheta\left[\begin{array}{cc} 0 \ 0 \\ 1 \ 1 \\ \end{array}\right] (0 , 0) \ 
 \pm \varTheta\left[\begin{array}{cc} 1 \ 1 \\ 1 \ 1 \\ \end{array}\right] (0 , 0) \Big)^2 
\nonumber\\
&&\hskip-15mm 
=\vartheta\left[\begin{array}{cc} 0 \ 0 \\ 1 \ 1 \\ \end{array}\right]  (0 , 0) \  
 \vartheta\left[\begin{array}{cc} 0 \ 0 \\ 0 \ 0 \\ \end{array}\right] (0 , 0) \ \
\pm\vartheta\left[\begin{array}{cc} 1 \ 1 \\ 1 \ 1 \\ \end{array}\right] (0 , 0) \ \
\vartheta\left[\begin{array}{cc} 1 \ 1 \\ 0 \ 0 \\ \end{array}\right] (0 , 0)  .
\label{2e41}
\end{eqnarray} 
In this way, from Eq.(\ref{2e36}), Eq.(\ref{2e38}), Eq.(\ref{2e39}), Eq.(\ref{2e41}), 
we have expressions 
$\varTheta\left[\begin{array}{cc} * \ * \\ * \ * \\ \end{array}\right] (0 ,0)$ by
$\vartheta\left[\begin{array}{cc} * \ * \\ * \ * \\ \end{array}\right] (0 ,0)$.
For six odd functions at $u=0$, $v=0$, 
$\varTheta\left[\begin{array}{cc} 0 \ 1 \\ 0 \ 1 \\ \end{array}\right] (0 ,0)$,\ 
$\varTheta\left[\begin{array}{cc} 1 \ 0 \\ 1 \ 0 \\ \end{array}\right] (0 ,0)$,\ 
$\varTheta\left[\begin{array}{cc} 0 \ 1 \\ 1 \ 1 \\ \end{array}\right] (0 ,0)$,\ 
$\varTheta\left[\begin{array}{cc} 1 \ 0 \\ 1 \ 1 \\ \end{array}\right] (0 ,0)$,\ 
$\varTheta\left[\begin{array}{cc} 1 \ 1 \\ 0 \ 1 \\ \end{array}\right] (0 ,0)$,\  
$\varTheta\left[\begin{array}{cc} 1 \ 1 \\ 1 \ 0 \\ \end{array}\right] (0 ,0)$
become zero. \\
Explicit expressions are given in Appendix D.\\
%
\subsection{Three additional formulae by connecting 
between A) , B), C) and D) sectors }
Up to this stage, we can express  
$\displaystyle{\varTheta\left[\begin{array}{cc} * \ * \\ * \ * \\ \end{array}\right] (2u , 2v)}$
by $\displaystyle{\vartheta\left[\begin{array}{cc} * \ * \\ * \ * \\ \end{array}\right] (u , v)}$ and
$\displaystyle{\vartheta\left[\begin{array}{cc} * \ * \\ * \ * \\ \end{array}\right] (0 , 0)}$.\\
For A) sector, we have three addition formulae of  hyperelliptic function by eliminating\\
$\displaystyle{\vartheta\left[\begin{array}{cc} 0 \ 0 \\ 0 \ 0 \\ \end{array}\right] (u_1-u_2 , v_1-v_2)}$.
Similarly for each B), C), D) sector, we have three addition formulae. 
Then we totally have twelve addition formulae of genus two hyperelliptic functions.
In order to obtain full addition formulae, we must connect between A), 
B), C) and D) sectors.\\
In order to obtain remaining three addition formulae, we must connect 
B), C), D) sectors with A) sector. 
\\
We use the term 
\begin{eqnarray}
&&\hskip-15mm \vartheta\left[\begin{array}{cc} 0 \ 1  \\ 0 \ 0 \\ \end{array}\right](u_1+u_2 ,v_1+v_2) \ 
\vartheta\left[\begin{array}{cc} 0 \ 0 \\ 0 \ 0 \\ \end{array}\right](u_1-u_2, v_1-v_2) ,
\nonumber
\end{eqnarray}
to connect B) sector with A) sector.
\\
We use the term
\begin{eqnarray}
&&\hskip-15mm \vartheta\left[\begin{array}{cc} 1 \ 0 \\ 0 \ 0 \\ \end{array}\right](u_1+u_2 ,v_1+v_2) \ 
\vartheta\left[\begin{array}{cc} 0 \ 0 \\ 0 \ 0 \\ \end{array}\right](u_1-u_2, v_1-v_2) ,
\nonumber
\end{eqnarray}
to connect C) sector with A) sector. 
\\
We use the term 
\begin{eqnarray}
&&\hskip-15mm \vartheta\left[\begin{array}{cc} 1 \ 1 \\ 0 \ 0 \\ \end{array}\right](u_1+u_2 ,v_1+v_2) \ 
\vartheta\left[\begin{array}{cc} 0 \ 0 \\ 0 \ 0 \\ \end{array}\right](u_1-u_2, v_1-v_2) ,
\nonumber
\end{eqnarray}
to connect D) sector with A) sector.
%
\subsection{Connection of D) sector with A) sector}
{\bf Addition formulae}\\
Using the term 
\begin{eqnarray}
&&\hskip-15mm \vartheta\left[\begin{array}{cc} 1 \ 1 \\ 0 \ 0 \\ \end{array}\right](u_1+u_2 ,v_1+v_2) \ 
\vartheta\left[\begin{array}{cc} 0 \ 0 \\ 0 \ 0 \\ \end{array}\right](u_1-u_2, v_1-v_2) , 
\nonumber
\end{eqnarray}
we can connect D) sector with A) sector.
We use the following duplication formula
\begin{eqnarray}
&&\hskip-15mm \vartheta\left[\begin{array}{cc} 1 \ 1 \\ 0 \ 0 \\ \end{array}\right](u_1+u_2 ,v_1+v_2) \ 
\vartheta\left[\begin{array}{cc} 0 \ 0 \\ 0 \ 0 \\ \end{array}\right](u_1-u_2, v_1-v_2) 
\nonumber\\
&&\hskip-15mm 
=\varTheta\left[\begin{array}{cc} \frac{1}{2} \ \frac{1}{2} \\ 0 \ 0 \\ \end{array}\right](2u_1 , 2v_1) \  
\varTheta\left[\begin{array}{cc} \frac{1}{2}  \ \frac{1}{2} \\ 0 \ 0 \\ \end{array}\right](2u_2 , 2v_2) 
+\varTheta\left[\begin{array}{cc} \frac{1}{2} \ -\frac{1}{2} \\ 0 \ 0 \\ \end{array}\right](2u_1 , 2v_1) \ 
\varTheta\left[\begin{array}{cc} \frac{1}{2} \ -\frac{1}{2} \\ 0 \ 0 \\ \end{array}\right](2u_2 , 2v_2) \
\nonumber\\
&&\hskip-15mm 
+\varTheta\left[\begin{array}{cc} -\frac{1}{2} \ \frac{1}{2} \\ 0 \ 0 \\ \end{array}\right](2u_1 , 2v_1) \ 
\varTheta\left[\begin{array}{cc} -\frac{1}{2} \ \frac{1}{2} \\ 0 \ 0 \\ \end{array}\right](2u_2 , 2v_2) 
+\varTheta\left[\begin{array}{cc} -\frac{1}{2} \ -\frac{1}{2} \\ 0 \ 0 \\ \end{array}\right](2u_1 , 2v_1) \
\varTheta\left[\begin{array}{cc} -\frac{1}{2} \ -\frac{1}{2} \\ 0 \ 0 \\ \end{array}\right](2u_2 , 2v_2)  . 
\nonumber\\
\label{2e42}
\end{eqnarray}
We take the ratio of Eq.(\ref{2e42}) with  Eq.(\ref{2e8}), we have the addition formula
of \\
$ \vartheta\left[\begin{array}{cc} 1 \ 1 \\ 0 \ 0 \\ \end{array}\right](u_1+u_2, v_1+v_2) /
\vartheta\left[\begin{array}{cc} 0 \ 0 \\ 0 \ 0 \\ \end{array}\right](u_1+u_2, _1+v_2)$.
Then we can connect addition formulae of A) sector with addition formulae of D) sector. 
Explicit expressions are given in appendix A by using Appendix B.\\
\\
{\bf Functional relations}\\
To obtain functional relations, we put $u_2=0$, $v_2=0$, $u_1=u$, $v_1=v$.
In the right-hand side of Eq.(\ref{2e42}),
$\varTheta\left[\begin{array}{cc} \frac{1}{2} \ \frac{1}{2} \\ 0 \ 0 \\ \end{array}\right](2u , 2v)$, 
$\varTheta\left[\begin{array}{cc} \frac{1}{2} \ -\frac{1}{2} \\ 0 \ 0 \\ \end{array}\right](2u , 2v)$, 
$\varTheta\left[\begin{array}{cc} -\frac{1}{2} \ \frac{1}{2} \\ 0 \ 0 \\ \end{array}\right](2u , 2v)$,\\ 
$\varTheta\left[\begin{array}{cc} -\frac{1}{2} \ -\frac{1}{2} \\ 0 \ 0 \\ \end{array}\right](2u , 2v)$, 
appears. Therefore, we must consider three more products, where the right hand side 
contains these four types.
Then we consider
\begin{eqnarray}
&&\hskip-15mm a)\ \vartheta\left[\begin{array}{cc} 1 \ 1  \\ 0 \ 0 \\ \end{array}\right](u ,v) \ 
\vartheta\left[\begin{array}{cc} 0 \ 0 \\ 0 \ 0 \\ \end{array}\right](u, v) ,
\label{2e43}\\
&&\hskip-15mm b)\ \vartheta\left[\begin{array}{cc} 0 \ 1  \\ 0 \ 0 \\ \end{array}\right](u ,v) \ 
\vartheta\left[\begin{array}{cc} 1 \ 0 \\ 0 \ 0 \\ \end{array}\right](u, v) ,
\label{2e44}\\
&&\hskip-15mm c)\ \vartheta\left[\begin{array}{cc} 1 \ 1  \\ 1 \ 0 \\ \end{array}\right](u ,v) \ 
\vartheta\left[\begin{array}{cc} 0 \ 0 \\ 1 \ 0 \\ \end{array}\right](u, v) ,
\label{2e45}\\
&&\hskip-15mm d)\ \vartheta\left[\begin{array}{cc} 1 \ 0  \\ 1 \ 0 \\ \end{array}\right](u ,v) \ 
\vartheta\left[\begin{array}{cc} 0 \ 1 \\ 1 \ 0 \\ \end{array}\right](u, v)  .
\label{2e46}
\end{eqnarray}
\\
For each term, we have 
\begin{eqnarray}
&&\hskip-15mm a)\ \vartheta\left[\begin{array}{cc} 1 \ 1 \\ 0 \ 0 \\ \end{array}\right](u ,v) \ 
\vartheta\left[\begin{array}{cc} 0 \ 0 \\ 0 \ 0 \\ \end{array}\right](u,v) 
\nonumber\\
&&\hskip-15mm 
=\Big( \varTheta\left[\begin{array}{cc} \frac{1}{2} \ \frac{1}{2} \\ 0 \ 0 \\ \end{array}\right](2u , 2v) 
+\varTheta\left[\begin{array}{cc} -\frac{1}{2} \ -\frac{1}{2} \\ 0 \ 0 \\ \end{array}\right](2u , 2v)\Big)
\varTheta\left[\begin{array}{cc} \frac{1}{2}  \ \frac{1}{2} \\ 0 \ 0 \\ \end{array}\right](0 , 0) 
\nonumber\\
&&\hskip-15mm 
+\Big( \varTheta\left[\begin{array}{cc} \frac{1}{2} \ -\frac{1}{2} \\ 0 \ 0 \\ \end{array}\right](2u , 2v) 
+\varTheta\left[\begin{array}{cc} -\frac{1}{2} \ \frac{1}{2} \\ 0 \ 0 \\ \end{array}\right](2u , 2v)\Big)
\varTheta\left[\begin{array}{cc} \frac{1}{2}  \ -\frac{1}{2} \\ 0 \ 0 \\ \end{array}\right](0 , 0) , 
\label{2e47}
\end{eqnarray}
\begin{eqnarray}
&&\hskip-15mm b)\ \vartheta\left[\begin{array}{cc} 0 \ 1 \\ 0 \ 0 \\ \end{array}\right](u ,v) \ 
\vartheta\left[\begin{array}{cc} 1 \ 0 \\ 0 \ 0 \\ \end{array}\right](u,v) 
\nonumber\\
&&\hskip-15mm 
=\Big( \varTheta\left[\begin{array}{cc} \frac{1}{2} \ -\frac{1}{2} \\ 0 \ 0 \\ \end{array}\right](2u , 2v) 
+\varTheta\left[\begin{array}{cc} -\frac{1}{2} \ \frac{1}{2} \\ 0 \ 0 \\ \end{array}\right](2u , 2v)\Big)
\varTheta\left[\begin{array}{cc} \frac{1}{2}  \ \frac{1}{2} \\ 0 \ 0 \\ \end{array}\right](0 , 0) 
\nonumber\\
&&\hskip-15mm 
+\Big( \varTheta\left[\begin{array}{cc} \frac{1}{2} \ \frac{1}{2} \\ 0 \ 0 \\ \end{array}\right](2u , 2v) 
+\varTheta\left[\begin{array}{cc} -\frac{1}{2} \ -\frac{1}{2} \\ 0 \ 0 \\ \end{array}\right](2u , 2v)\Big)
\varTheta\left[\begin{array}{cc} \frac{1}{2}  \ -\frac{1}{2} \\ 0 \ 0 \\ \end{array}\right](0 , 0) , 
\label{2e48}
\end{eqnarray}
\begin{eqnarray}
&&\hskip-15mm c)\ \vartheta\left[\begin{array}{cc} 1 \ 1 \\ 1 \ 0 \\ \end{array}\right](u ,v) \ 
\vartheta\left[\begin{array}{cc} 0 \ 0 \\ 1 \ 0 \\ \end{array}\right](u,v) 
\nonumber\\
&&\hskip-15mm 
=\sqrt{-1}\Big( \varTheta\left[\begin{array}{cc} \frac{1}{2} \ \frac{1}{2} \\ 0 \ 0 \\ 
\end{array}\right](2u , 2v) 
-\varTheta\left[\begin{array}{cc} -\frac{1}{2} \ -\frac{1}{2} \\ 0 \ 0 \\ \end{array}\right](2u , 2v)\Big)
\varTheta\left[\begin{array}{cc} \frac{1}{2}  \ \frac{1}{2} \\ 0 \ 0 \\ \end{array}\right](0 , 0) 
\nonumber\\
&&\hskip-15mm 
+\sqrt{-1}\Big( \varTheta\left[\begin{array}{cc} \frac{1}{2} \ -\frac{1}{2} \\ 0 \ 0 \\ 
\end{array}\right](2u , 2v) 
-\varTheta\left[\begin{array}{cc} -\frac{1}{2} \ \frac{1}{2} \\ 0 \ 0 \\ \end{array}\right](2u , 2v)\Big)
\varTheta\left[\begin{array}{cc} \frac{1}{2}  \ -\frac{1}{2} \\ 0 \ 0 \\ \end{array}\right](0 , 0) ,
\label{2e49}
\end{eqnarray}
\begin{eqnarray}
&&\hskip-15mm d)\ \vartheta\left[\begin{array}{cc} 1 \ 0 \\ 1 \ 0 \\ \end{array}\right](u ,v) \ 
\vartheta\left[\begin{array}{cc} 0 \ 1 \\ 1 \ 0 \\ \end{array}\right](u,v) 
\nonumber\\
&&\hskip-15mm 
=\sqrt{-1}\Big( \varTheta\left[\begin{array}{cc} \frac{1}{2} \ \frac{1}{2} \\ 0 \ 0 \\ 
\end{array}\right](2u , 2v) 
-\varTheta\left[\begin{array}{cc} -\frac{1}{2} \ -\frac{1}{2} \\ 0 \ 0 \\ \end{array}\right](2u , 2v)\Big)
\varTheta\left[\begin{array}{cc} \frac{1}{2}  \ -\frac{1}{2} \\ 0 \ 0 \\ \end{array}\right](0 , 0) 
\nonumber\\
&&\hskip-15mm 
+\sqrt{-1}\Big( \varTheta\left[\begin{array}{cc} \frac{1}{2} \ -\frac{1}{2} \\ 0 \ 0 \\ 
\end{array}\right](2u , 2v) 
-\varTheta\left[\begin{array}{cc} -\frac{1}{2} \ \frac{1}{2} \\ 0 \ 0 \\ \end{array}\right](2u , 2v)\Big)
\varTheta\left[\begin{array}{cc} \frac{1}{2}  \ \frac{1}{2} \\ 0 \ 0 \\ \end{array}\right](0 , 0) .
\label{2e50}
\end{eqnarray}
where we use 
$\varTheta\left[\begin{array}{cc} \frac{1}{2}  \ \frac{1}{2} \\ 0 \ 0 \\ \end{array}\right](0 , 0)
=\varTheta\left[\begin{array}{cc} -\frac{1}{2}  \ -\frac{1}{2} \\ 0 \ 0 \\ \end{array}\right](0 , 0)$,
$\varTheta\left[\begin{array}{cc} \frac{1}{2}  \ -\frac{1}{2} \\ 0 \ 0 \\ \end{array}\right](0 , 0)
=\varTheta\left[\begin{array}{cc} -\frac{1}{2}  \ \frac{1}{2} \\ 0 \ 0 \\ \end{array}\right](0 , 0)$.\\
Using Eq.(\ref{2e47}), Eq.(\ref{2e48}), Eq.(\ref{2e49}), Eq.(\ref{2e50}), we can express
$\varTheta\left[\begin{array}{cc} \frac{1}{2} \ \frac{1}{2} \\ 0 \ 0 \\ \end{array}\right](2u , 2v)$, \\
$\varTheta\left[\begin{array}{cc} \frac{1}{2} \ -\frac{1}{2} \\ 0 \ 0 \\ \end{array}\right](2u , 2v)$, 
$\varTheta\left[\begin{array}{cc} -\frac{1}{2} \ \frac{1}{2} \\ 0 \ 0 \\ \end{array}\right](2u , 2v)$, 
$\varTheta\left[\begin{array}{cc} -\frac{1}{2} \ -\frac{1}{2} \\ 0 \ 0 \\ \end{array}\right](2u , 2v)$ 
by $\vartheta\left[\begin{array}{cc} * \ * \\ * \ * \\ \end{array}\right](u , v)$.\\
Explicit expressions are given in Appendix C.\\
\\
{\bf Constants relations}\\
By putting $u=0$, $v=0$ in Eq.(\ref{2e47}) and Eq.(\ref{2e48}), we have constants relations
\begin{eqnarray}
&&\hskip-25mm a')\ \vartheta\left[\begin{array}{cc} 1 \ 1 \\ 0 \ 0 \\ \end{array}\right](0 ,0) \ 
\vartheta\left[\begin{array}{cc} 0 \ 0 \\ 0 \ 0 \\ \end{array}\right](0,0) 
=2\Big( \varTheta\left[\begin{array}{cc} \frac{1}{2} \ \frac{1}{2} \\ 0 \ 0 \\ \end{array}\right]^2 (0 , 0) 
+\varTheta\left[\begin{array}{cc} \frac{1}{2} \ -\frac{1}{2} \\ 0 \ 0 \\ \end{array}\right]^2 (0 , 0)\Big) ,
\label{2e51}\\
%
&&\hskip-25mm b')\ \vartheta\left[\begin{array}{cc} 0 \ 1 \\ 0 \ 0 \\ \end{array}\right](0 ,0) \ 
\vartheta\left[\begin{array}{cc} 1 \ 0 \\ 0 \ 0 \\ \end{array}\right](0,0) 
=4 \varTheta\left[\begin{array}{cc} \frac{1}{2} \ \frac{1}{2} \\ 0 \ 0 \\ \end{array}\right] (0 , 0) \ 
\varTheta\left[\begin{array}{cc} \frac{1}{2} \ -\frac{1}{2} \\ 0 \ 0 \\ \end{array}\right] (0 , 0)  .
\label{2e52}
\end{eqnarray}
Combining Eq.(\ref{2e51}) and Eq.(\ref{2e52}), we obtain 
\begin{eqnarray}
&&\hskip-15mm \vartheta\left[\begin{array}{cc} 0 \ 1 \\ 0 \ 0 \\ \end{array}\right](0 ,0) \ 
\vartheta\left[\begin{array}{cc} 1 \ 0 \\ 0 \ 0 \\ \end{array}\right](0,0) 
\pm \vartheta\left[\begin{array}{cc} 0 \ 1 \\ 0 \ 0 \\ \end{array}\right](0 ,0) \ 
\vartheta\left[\begin{array}{cc} 1 \ 0 \\ 0 \ 0 \\ \end{array}\right](0,0) 
\nonumber\\
&&\hskip-15mm 
=2 \Big( \varTheta\left[\begin{array}{cc} \frac{1}{2} \ \frac{1}{2} \\ 0 \ 0 \\ \end{array}\right] (0 , 0)  
\pm \varTheta\left[\begin{array}{cc} \frac{1}{2} \ -\frac{1}{2} \\ 0 \ 0 \\ \end{array}\right] (0 , 0)\Big)^2  .
\label{2e53}
\end{eqnarray}
Then
$\varTheta\left[\begin{array}{cc} \frac{1}{2} \ \frac{1}{2} \\ 0 \ 0 \\ \end{array}\right] (0 , 0)
=\varTheta\left[\begin{array}{cc} -\frac{1}{2} \ -\frac{1}{2} \\ 0 \ 0 \\ \end{array}\right] (0 , 0)$,
$\varTheta\left[\begin{array}{cc} \frac{1}{2} \ -\frac{1}{2} \\ 0 \ 0 \\ \end{array}\right] (0 , 0)
=\varTheta\left[\begin{array}{cc} -\frac{1}{2} \ \frac{1}{2} \\ 0 \ 0 \\ \end{array}\right] (0 , 0)$,    
are given by $\vartheta\left[\begin{array}{cc} * \ * \\ * \ * \\ \end{array}\right] (0 , 0)$.
Explicit expressions are given by Appendix D.
%
\subsection{Connection of B) sector with A) sector}
{\bf Addition formula}\\
Using the term
\begin{eqnarray}
&&\hskip-15mm \vartheta\left[\begin{array}{cc} 0 \ 1  \\ 0 \ 0 \\ \end{array}\right](u_1+u_2 ,v_1+v_2) \ 
\vartheta\left[\begin{array}{cc} 0 \ 0 \\ 0 \ 0 \\ \end{array}\right](u_1-u_2, v_1-v_2)  ,
\nonumber
\end{eqnarray}
we can connect B) sector with A) sector. We use the following duplication formula 
\begin{eqnarray}
&&\hskip-15mm \vartheta\left[\begin{array}{cc} 0 \ 1 \\ 0 \ 0 \\ \end{array}\right](u_1+u_2 ,v_1+v_2) \ 
\vartheta\left[\begin{array}{cc} 0 \ 0 \\ 0 \ 0 \\ \end{array}\right](u_1-u_2, v_1-v_2) 
\nonumber\\
&&\hskip-15mm 
=\varTheta\left[\begin{array}{cc} 0 \ \frac{1}{2} \\ 0 \ 0 \\ \end{array}\right](2u_1 , 2v_1) \ \ 
\varTheta\left[\begin{array}{cc} 0  \ \frac{1}{2} \\ 0 \ 0 \\ \end{array}\right](2u_2 , 2v_2) 
+\varTheta\left[\begin{array}{cc} 0 \ -\frac{1}{2} \\ 0 \ 0 \\ \end{array}\right](2u_1 , 2v_1) \ \
\varTheta\left[\begin{array}{cc} 0 \ -\frac{1}{2} \\ 0 \ 0 \\ \end{array}\right](2u_2 , 2v_2) \
\nonumber\\
&&\hskip-15mm 
+\varTheta\left[\begin{array}{cc} 1 \ \frac{1}{2} \\ 0 \ 0 \\ \end{array}\right](2u_1 , 2v_1) \ \
\varTheta\left[\begin{array}{cc} 1 \ \frac{1}{2} \\ 0 \ 0 \\ \end{array}\right](2u_2 , 2v_2) 
+\varTheta\left[\begin{array}{cc} 1 \ -\frac{1}{2} \\ 0 \ 0 \\ \end{array}\right](2u_1 , 2v_1) \
\varTheta\left[\begin{array}{cc} 1 \ -\frac{1}{2} \\ 0 \ 0 \\ \end{array}\right](2u_2 , 2v_2) . 
\nonumber\\
\label{2e54}
\end{eqnarray}
We take the ratio of Eq.(\ref{2e54}) with  Eq.(\ref{2e8}), we have the addition formula
of \\
$ \vartheta\left[\begin{array}{cc} 0 \ 1 \\ 0 \ 0 \\ \end{array}\right](u_1+u_2, v_1+v_2) /
\vartheta\left[\begin{array}{cc} 0 \ 0 \\ 0 \ 0 \\ \end{array}\right](u_1+u_2, _1+v_2)$. 
Then we can connect addition formulae of A) sector with addition formulae of 
D) sector.  Explicit expressions are given 
in appendix A by using Appendix B.\\
\\
{\bf Functional relations}\\
To obtain the function relation, we put $u_2=0$, $v_2=0$, $u_1=u$, $v_1=v$.
In the right-hand side,
$\varTheta\left[\begin{array}{cc} 0 \ \frac{1}{2} \\ 0 \ 0 \\ \end{array}\right](2u , 2v)$, 
$\varTheta\left[\begin{array}{cc} 0 \ -\frac{1}{2} \\ 0 \ 0 \\ \end{array}\right](2u , 2v)$, 
$\varTheta\left[\begin{array}{cc} 1 \ \frac{1}{2} \\ 0 \ 0 \\ \end{array}\right](2u , 2v)$,
$\varTheta\left[\begin{array}{cc} 1 \ -\frac{1}{2} \\ 0 \ 0 \\ \end{array}\right](2u , 2v)$ \\
appears. Therefore, we must consider three more products, where the right hand side 
contains these four types.
Then we consider
\begin{eqnarray}
&&\hskip-15mm a)\ \vartheta\left[\begin{array}{cc} 0 \ 1  \\ 0 \ 0 \\ \end{array}\right](u ,v) \ 
\vartheta\left[\begin{array}{cc} 0 \ 0 \\ 0 \ 0 \\ \end{array}\right](u, v) ,
\label{2e55}\\
&&\hskip-15mm b)\ \vartheta\left[\begin{array}{cc} 0 \ 1  \\ 1 \ 0 \\ \end{array}\right](u ,v) \ 
\vartheta\left[\begin{array}{cc} 0 \ 0 \\ 1 \ 0 \\ \end{array}\right](u, v) ,
\label{2e56}\\
&&\hskip-15mm c)\ \vartheta\left[\begin{array}{cc} 0 \ 1  \\ 0 \ 1 \\ \end{array}\right](u ,v) \ 
\vartheta\left[\begin{array}{cc} 0 \ 0 \\ 0 \ 1 \\ \end{array}\right](u, v) ,
\label{2e57}\\
&&\hskip-15mm d)\ \vartheta\left[\begin{array}{cc} 0 \ 1  \\ 1 \ 1 \\ \end{array}\right](u ,v) \ 
\vartheta\left[\begin{array}{cc} 0 \ 0 \\ 1 \ 1 \\ \end{array}\right](u, v) .
\label{2e58}
\end{eqnarray}
For each term, we have 
\begin{eqnarray}
&&\hskip-15mm a)\ \vartheta\left[\begin{array}{cc} 0 \ 1 \\ 0 \ 0 \\ \end{array}\right](u ,v) \ 
\vartheta\left[\begin{array}{cc} 0 \ 0 \\ 0 \ 0 \\ \end{array}\right](0,0) 
\nonumber\\
&&\hskip-15mm 
=\Big( \varTheta\left[\begin{array}{cc} 0 \ \frac{1}{2} \\ 0 \ 0 \\ \end{array}\right](2u , 2v) 
+\varTheta\left[\begin{array}{cc} 0  \  -\frac{1}{2} \\ 0 \ 0 \\ \end{array}\right](2u , 2v)\Big)
\varTheta\left[\begin{array}{cc} 0 \ \frac{1}{2} \\ 0 \ 0 \\ \end{array}\right](0 , 0) 
\nonumber\\
&&\hskip-15mm 
+\Big( \varTheta\left[\begin{array}{cc} 1 \ \frac{1}{2} \\ 0 \ 0 \\ \end{array}\right](2u , 2v) 
+\varTheta\left[\begin{array}{cc} 1 \ -\frac{1}{2} \\ 0 \ 0 \\ \end{array}\right](2u_1 , 2v_1)\Big)
\varTheta\left[\begin{array}{cc} 1  \ \frac{1}{2} \\ 0 \ 0 \\ \end{array}\right](0 , 0) ,
\label{2e59}
\end{eqnarray}
\begin{eqnarray}
&&\hskip-15mm b)\ \vartheta\left[\begin{array}{cc} 0 \ 1 \\ 1 \ 0 \\ \end{array}\right](u ,v) \ 
\vartheta\left[\begin{array}{cc} 0 \ 0 \\ 1 \ 0 \\ \end{array}\right](0,0) 
\nonumber\\
&&\hskip-15mm 
=\Big( \varTheta\left[\begin{array}{cc} 0 \ \frac{1}{2} \\ 0 \ 0 \\ \end{array}\right](2u , 2v) 
+\varTheta\left[\begin{array}{cc} 0  \  -\frac{1}{2} \\ 0 \ 0 \\ \end{array}\right](2u , 2v)\Big)
\varTheta\left[\begin{array}{cc} 0 \ \frac{1}{2} \\ 0 \ 0 \\ \end{array}\right](0 , 0) 
\nonumber\\
&&\hskip-15mm 
-\Big( \varTheta\left[\begin{array}{cc} 1 \ \frac{1}{2} \\ 0 \ 0 \\ \end{array}\right](2u_1 , 2v_1) 
+\varTheta\left[\begin{array}{cc} 1 \ -\frac{1}{2} \\ 0 \ 0 \\ \end{array}\right](2u_1 , 2v_1)\Big)
\varTheta\left[\begin{array}{cc} 1  \ \frac{1}{2} \\ 0 \ 0 \\ \end{array}\right](0 , 0) ,
\label{2e60}
\end{eqnarray}
\begin{eqnarray}
&&\hskip-15mm c)\ \vartheta\left[\begin{array}{cc} 0 \ 1 \\ 0 \ 1 \\ \end{array}\right](u ,v) \ 
\vartheta\left[\begin{array}{cc} 0 \ 0 \\ 0 \ 1 \\ \end{array}\right](u,v) 
\nonumber\\
&&\hskip-15mm 
=\sqrt{-1}\Big( \varTheta\left[\begin{array}{cc} 0 \ \frac{1}{2} \\ 0 \ 0 \\ 
\end{array}\right](2u , 2v) 
-\varTheta\left[\begin{array}{cc} 0 \ -\frac{1}{2} \\ 0 \ 0 \\ \end{array}\right](2u , 2v)\Big)
\varTheta\left[\begin{array}{cc} 0  \ \frac{1}{2} \\ 0 \ 0 \\ \end{array}\right](0 , 0) 
\nonumber\\
&&\hskip-15mm 
+\sqrt{-1}\Big( \varTheta\left[\begin{array}{cc} 1 \ \frac{1}{2} \\ 0 \ 0 \\ 
\end{array}\right](2u , 2v) 
-\varTheta\left[\begin{array}{cc} 1 \ -\frac{1}{2} \\ 0 \ 0 \\ \end{array}\right](2u , 2v)\Big)
\varTheta\left[\begin{array}{cc} 1  \ \frac{1}{2} \\ 0 \ 0 \\ \end{array}\right](0 , 0) ,
\label{2e61}
\end{eqnarray}
\begin{eqnarray}
&&\hskip-15mm d)\ \vartheta\left[\begin{array}{cc} 0 \ 1 \\ 1 \ 1 \\ \end{array}\right](u ,v) \ 
\vartheta\left[\begin{array}{cc} 0 \ 0 \\ 1 \ 1 \\ \end{array}\right](u,v) 
\nonumber\\
&&\hskip-15mm 
=\sqrt{-1}\Big( \varTheta\left[\begin{array}{cc} 0 \ \frac{1}{2} \\ 0 \ 0 \\ 
\end{array}\right](2u , 2v) 
-\varTheta\left[\begin{array}{cc} 0 \ -\frac{1}{2} \\ 0 \ 0 \\ \end{array}\right](2u , 2v)\Big)
\varTheta\left[\begin{array}{cc} 0  \ \frac{1}{2} \\ 0 \ 0 \\ \end{array}\right](0 , 0) 
\nonumber\\
&&\hskip-15mm 
-\sqrt{-1}\Big( \varTheta\left[\begin{array}{cc} 1 \ \frac{1}{2} \\ 0 \ 0 \\ 
\end{array}\right](2u , 2v) 
-\varTheta\left[\begin{array}{cc} 1 \ -\frac{1}{2} \\ 0 \ 0 \\ \end{array}\right](2u , 2v)\Big)
\varTheta\left[\begin{array}{cc} 1  \ \frac{1}{2} \\ 0 \ 0 \\ \end{array}\right](0 , 0) ,
\label{2e62}
\end{eqnarray}
where we use 
$\varTheta\left[\begin{array}{cc} 0  \ \frac{1}{2} \\ 0 \ 0 \\ \end{array}\right](0 , 0)
=\varTheta\left[\begin{array}{cc} 0  \ -\frac{1}{2} \\ 0 \ 0 \\ \end{array}\right](0 , 0)$,
$\varTheta\left[\begin{array}{cc} 1  \ \frac{1}{2} \\ 0 \ 0 \\ \end{array}\right](0 , 0)
=\varTheta\left[\begin{array}{cc} 1  \ -\frac{1}{2} \\ 0 \ 0 \\ \end{array}\right](0 , 0)$,
Using Eq.(\ref{2e59}), Eq.(\ref{2e60}), Eq.(\ref{2e61}), Eq.(\ref{2e62}), we can express
$\varTheta\left[\begin{array}{cc} 0 \ \frac{1}{2} \\ 0 \ 0 \\ \end{array}\right](2u , 2v)$, \\ 
$\varTheta\left[\begin{array}{cc} 0 \ -\frac{1}{2} \\ 0 \ 0 \\ \end{array}\right](2u , 2v)$, 
$\varTheta\left[\begin{array}{cc} 1 \ \frac{1}{2} \\ 0 \ 0 \\ \end{array}\right](2u , 2v)$, 
$\varTheta\left[\begin{array}{cc} 1 \ -\frac{1}{2} \\ 0 \ 0 \\ \end{array}\right](2u , 2v)$ 
by $\vartheta\left[\begin{array}{cc} * \ * \\ * \ * \\ \end{array}\right](u , v)$.\\
Explicit expressions are given in Appendix C.\\
\\
{\bf Constants relations} \\
By putting $u=0$, $v=0$ in Eq.(\ref{2e59}) and Eq.(\ref{2e60}), we have 
\begin{eqnarray}
&&\hskip-25mm a')\ \vartheta\left[\begin{array}{cc} 0 \ 1 \\ 0 \ 0 \\ \end{array}\right](0 ,0) \ 
\vartheta\left[\begin{array}{cc} 0 \ 0 \\ 0 \ 0 \\ \end{array}\right](0,0) 
=2\Big( \varTheta\left[\begin{array}{cc} 0 \ \frac{1}{2} \\ 0 \ 0 \\ \end{array}\right]^2 (0 , 0) 
+\varTheta\left[\begin{array}{cc} 1 \ \frac{1}{2} \\ 0 \ 0 \\ \end{array}\right]^2 (0 , 0)\Big) ,
\label{2e63}
\end{eqnarray}
\begin{eqnarray}
&&\hskip-25mm b')\ \vartheta\left[\begin{array}{cc} 0 \ 1 \\ 1 \ 0 \\ \end{array}\right](0 ,0) \ 
\vartheta\left[\begin{array}{cc} 0 \ 0 \\ 1 \ 0 \\ \end{array}\right](0,0) 
=2\Big( \varTheta\left[\begin{array}{cc} 0 \ \frac{1}{2} \\ 0 \ 0 \\ \end{array}\right]^2 (0 , 0) 
-\varTheta\left[\begin{array}{cc} 1 \ \frac{1}{2} \\ 0 \ 0 \\ \end{array}\right]^2 (0 , 0)\Big)  .
\label{2e64}
\end{eqnarray}
Then
$\varTheta\left[\begin{array}{cc} 0 \ \frac{1}{2} \\ 0 \ 0 \\ \end{array}\right] (0 , 0)
=\varTheta\left[\begin{array}{cc} 0 \ -\frac{1}{2} \\ 0 \ 0 \\ \end{array}\right] (0 , 0)$,
$\varTheta\left[\begin{array}{cc} 1 \ -\frac{1}{2} \\ 0 \ 0 \\ \end{array}\right] (0 , 0)
=\varTheta\left[\begin{array}{cc} 1 \ \frac{1}{2} \\ 0 \ 0 \\ \end{array}\right] (0 , 0)$    
is given by $\vartheta\left[\begin{array}{cc} a \ c \\ b \ d \\ \end{array}\right] (0 , 0)$.
Explicit expressions are given by Appendix D.
%
\subsection{Connection of C) sector with A) sector}
{\bf Addition formula}\\
We use the term 
\begin{eqnarray}
&&\hskip-15mm \vartheta\left[\begin{array}{cc} 1 \ 0 \\ 0 \ 0 \\ \end{array}\right](u_1+u_2 ,v_1+v_2) \ 
\vartheta\left[\begin{array}{cc} 0 \ 0 \\ 0 \ 0 \\ \end{array}\right](u_1-u_2, v_1-v_2) ,
\nonumber
\end{eqnarray}
to connect C) sector with A) sector. Formulae of the connection of C) sector with A) sector are 
given by the left-right symmetry from formulae of the connection of B) sector with A) sector.\\
From the left-right symmetry, we have
\begin{eqnarray}
&&\hskip-15mm \vartheta\left[\begin{array}{cc} 1 \ 0 \\ 0 \ 0 \\ \end{array}\right](u_1+u_2 ,v_1+v_2) \ 
\vartheta\left[\begin{array}{cc} 0 \ 0 \\ 0 \ 0 \\ \end{array}\right](u_1-u_2, v_1-v_2) 
\nonumber\\
&&\hskip-15mm 
=\varTheta\left[\begin{array}{cc} \frac{1}{2} \ 0 \\ 0 \ 0 \\ \end{array}\right](2u_1 , 2v_1) \ \ 
\varTheta\left[\begin{array}{cc}   \frac{1}{2} \ 0  \\ 0 \ 0 \\ \end{array}\right](2u_2 , 2v_2) 
+\varTheta\left[\begin{array}{cc} -\frac{1}{2} \ 0  \\ 0 \ 0 \\ \end{array}\right](2u_1 , 2v_1) \ \
\varTheta\left[\begin{array}{cc} -\frac{1}{2}  \ 0  \\ 0 \ 0 \\ \end{array}\right](2u_2 , 2v_2) \
\nonumber\\
&&\hskip-15mm 
+\varTheta\left[\begin{array}{cc} \frac{1}{2} \ 1 \\ 0 \ 0 \\ \end{array}\right](2u_1 , 2v_1) \ \
\varTheta\left[\begin{array}{cc}  \frac{1}{2} \ 1 \\ 0 \ 0 \\ \end{array}\right](2u_2 , 2v_2) 
+\varTheta\left[\begin{array}{cc} -\frac{1}{2} \ 1  \\ 0 \ 0 \\ \end{array}\right](2u_1 , 2v_1) \
\varTheta\left[\begin{array}{cc} -\frac{1}{2}  \ 1 \\ 0 \ 0 \\ \end{array}\right](2u_2 , 2v_2) .
\nonumber\\
\label{2e65}
\end{eqnarray}
\\
We take the ratio of Eq.(\ref{2e65}) with  Eq.(\ref{2e8}), we have the addition formula
of \\
$ \vartheta\left[\begin{array}{cc} 1 \ 0 \\ 0 \ 0 \\ \end{array}\right](u_1+u_2, v_1+v_2) /
\vartheta\left[\begin{array}{cc} 0 \ 0 \\ 0 \ 0 \\ \end{array}\right](u_1+u_2, _1+v_2)$. Then 
we can connect addition formulae of A) sector with addition formulae of C) sector. Explicit 
expressions are given 
in appendix A by using Appendix B.\\
\\
{\bf Functional relations}\\
To obtain the function relation, we put $u_2=0$, $v_2=0$, $u_1=u$, $v_1=v$.
In the right-hand side,
$\varTheta\left[\begin{array}{cc} \frac{1}{2} \ 0 \\ 0 \ 0 \\ \end{array}\right](2u , 2v)$, 
$\varTheta\left[\begin{array}{cc} -\frac{1}{2} \ 0  \\ 0 \ 0 \\ \end{array}\right](2u , 2v)$, 
$\varTheta\left[\begin{array}{cc} \frac{1}{2} \ 1 \\ 0 \ 0 \\ \end{array}\right](2u , 2v)$, 
$\varTheta\left[\begin{array}{cc} -\frac{1}{2} \ 1  \\ 0 \ 0 \\ \end{array}\right](2u , 2v)$, 
appears. Therefore, we must consider three more products, where the right hand side 
contains these four types.
Then we consider
\begin{eqnarray}
&&\hskip-15mm a)\ \vartheta\left[\begin{array}{cc} 1 \ 0  \\ 0 \ 0 \\ \end{array}\right](u ,v) \ 
\vartheta\left[\begin{array}{cc} 0 \ 0 \\ 0 \ 0 \\ \end{array}\right](u, v) ,
\label{2e66}\\
&&\hskip-15mm b)\ \vartheta\left[\begin{array}{cc} 1 \ 0  \\ 0 \ 1 \\ \end{array}\right](u ,v) \ 
\vartheta\left[\begin{array}{cc} 0 \ 0 \\ 0 \ 1 \\ \end{array}\right](u, v) ,
\label{2e67}\\
&&\hskip-15mm c)\ \vartheta\left[\begin{array}{cc} 1 \ 0  \\ 1 \ 0 \\ \end{array}\right](u ,v) \ 
\vartheta\left[\begin{array}{cc} 0 \ 0 \\ 1 \ 0 \\ \end{array}\right](u, v) ,
\label{2e68}\\
&&\hskip-15mm d)\ \vartheta\left[\begin{array}{cc} 1 \ 0  \\ 1 \ 1 \\ \end{array}\right](u ,v) \ 
\vartheta\left[\begin{array}{cc} 0 \ 0 \\ 1 \ 1 \\ \end{array}\right](u, v) .
\label{2e69}
\end{eqnarray}
For each term, we have 
\begin{eqnarray}
&&\hskip-15mm a)\ \vartheta\left[\begin{array}{cc} 1 \ 0 \\ 0 \ 0 \\ \end{array}\right](u ,v) \ 
\vartheta\left[\begin{array}{cc} 0 \ 0 \\ 0 \ 0 \\ \end{array}\right](0,0) 
\nonumber\\
&&\hskip-15mm 
=\Big( \varTheta\left[\begin{array}{cc} \frac{1}{2} \ 0 \\ 0 \ 0 \\ \end{array}\right](2u , 2v) 
+\varTheta\left[\begin{array}{cc} -\frac{1}{2} \ 0  \\ 0 \ 0 \\ \end{array}\right](2u , 2v)\Big)
\varTheta\left[\begin{array}{cc} \frac{1}{2}  \ 0 \\ 0 \ 0 \\ \end{array}\right](0 , 0) 
\nonumber\\
&&\hskip-15mm 
+\Big( \varTheta\left[\begin{array}{cc} \frac{1}{2} \ 1 \\ 0 \ 0 \\ \end{array}\right](2u , 2v) 
+\varTheta\left[\begin{array}{cc} -\frac{1}{2} \ 1  \\ 0 \ 0 \\ \end{array}\right](2u_1 , 2v_1)\Big)
\varTheta\left[\begin{array}{cc} \frac{1}{2} \ 1  \\ 0 \ 0 \\ \end{array}%
\right](0 , 0) ,
\label{2e73}
\end{eqnarray}
where we use 
$\varTheta\left[%
\right](0 , 0)$.
\\
Using Eq.(\ref{2e70}), Eq.(\ref{2e71}), Eq.(\ref{2e72}), Eq.(\ref{2e73}), we can express
$\varTheta\left[%
\right](u , v)$.\\
Explicit expressions are given in Appendix C.\\
\\
{\bf Constants relations} \\
By putting $u=0$, $v=0$ in Eq.(\ref{2e70}) and Eq.(\ref{2e71}), we have 
\begin{eqnarray}
&&\hskip-25mm a')\ \vartheta\left[%
\right] (0 , 0)$     
is given by $\vartheta\left[%
\right] (0 , 0)$.
Explicit expressions are given in Appendix D.
%
%
%
\vspace{5mm}
\section{Summary and Discussion}
\setcounter{equation}{0}
Here have given full algebraic addition formulae of genus two hyperelliptic functions
by the duplication method. 
Compared with the previous full algebraic addition formulae\cite{Shigemoto4} 
according to Kossak's method, this full algebraic addition formulae
according to the duplication method give some hint to the Lie group structure
of the addition formulae of genus two hyperelliptic functions.
In this derivation, the Riemann matrix naturally appears in Eq.(\ref{2e28}), Eq.(\ref{2e29}).
This Riemann matrix is connected with the representation of the 
special $SO(3)$ rotation expressed by the quaternion\cite{Hudson}.
This $SO(3)$ structure of the addition formulae of genus two hyper elliptic theta 
functions are essential to find the parametrization of the Kummer 
surface\cite{Gopel1,Gopel2,Borchardt,Cayley,Hudson}.\\
We have shown that there exists $SO(3)$ structure 
of the addition formulae of elliptic functions\cite{Shigemoto1,Shigemoto2}. 
(See Appendix E) Then we expect that there will exist some Lie group 
structure even for addition formulae of genus two hyperelliptic functions.

%
%
%
%
%
%
%
\newpage
%
%
%
%

\newpage
\appendix
\section{Addition formulae of hyperelliptic functions}
\setcounter{equation}{0}
\begin{eqnarray}
&&\hskip-15mm F[\begin{array}{cc} 0 \ 0 \\ 0 \ 1 \\ \end{array}](u_1+u_2 ,v_1+v_2)
=\frac{ \vartheta\left[\begin{array}{cc} 0 \ 0 \\ 0 \ 1 \\ \end{array}\right](u_1+u_2 ,v_1+v_2)} 
{\vartheta\left[\begin{array}{cc} 0 \ 0 \\ 0 \ 0 \\ \end{array}\right](u_1+u_2, v_1+v_2)} 
=\frac{ G[\begin{array}{cc} 0 \ 0 \\ 0 \ 1 \\ \end{array}] }
{G[\begin{array}{cc} 0 \ 0 \\ 0 \ 0 \\ \end{array}]}      ,
\label{A1}
\end{eqnarray}
\begin{eqnarray}
&&\hskip-15mm F[\begin{array}{cc} 0 \ 0 \\ 1 \ 0 \\ \end{array}](u_1+u_2 ,v_1+v_2)
=\frac{ \vartheta\left[\begin{array}{cc} 0 \ 0 \\ 1 \ 0 \\ \end{array}\right](u_1+u_2 ,v_1+v_2)} 
{\vartheta\left[\begin{array}{cc} 0 \ 0 \\ 0 \ 0 \\ \end{array}\right](u_1+u_2, v_1+v_2)} 
=\frac{ G[\begin{array}{cc} 0 \ 0 \\ 1 \ 0 \\ \end{array}] }
{G[\begin{array}{cc} 0 \ 0 \\ 0 \ 0 \\ \end{array}]}    ,
\label{A2}
\end{eqnarray}
\begin{eqnarray}
&&\hskip-15mm F[\begin{array}{cc} 0 \ 0 \\ 1 \ 1 \\ \end{array}](u_1+u_2 ,v_1+v_2)
=\frac{ \vartheta\left[\begin{array}{cc} 0 \ 0 \\ 1 \ 1 \\ \end{array}\right](u_1+u_2 ,v_1+v_2)} 
{\vartheta\left[\begin{array}{cc} 0 \ 0 \\ 0 \ 0 \\ \end{array}\right](u_1+u_2, v_1+v_2)} 
=\frac{ G[\begin{array}{cc} 0 \ 0 \\ 1 \ 1 \\ \end{array}] }
{G[\begin{array}{cc} 0 \ 0 \\ 0 \ 0 \\ \end{array}]}    ,
\label{A3}
\end{eqnarray}
\begin{eqnarray}
&&\hskip-15mm F[\begin{array}{cc} 0 \ 1 \\ 0 \ 0 \\ \end{array}](u_1+u_2 ,v_1+v_2)
=\frac{ \vartheta\left[\begin{array}{cc} 0 \ 1 \\ 0 \ 0 \\ \end{array}\right](u_1+u_2 ,v_1+v_2)} 
{\vartheta\left[\begin{array}{cc} 0 \ 0 \\ 0 \ 0 \\ \end{array}\right](u_1+u_2, v_1+v_2)} 
=\frac{ G_1[\begin{array}{cc} 0 \ 1 \\ 0 \ 0 \\ \end{array}] }
{G[\begin{array}{cc} 0 \ 0 \\ 0 \ 0 \\ \end{array}]}    ,
\label{A4}
\end{eqnarray}
\begin{eqnarray}
&&\hskip-15mm F[\begin{array}{cc} 0 \ 1 \\ 0 \ 1 \\ \end{array}](u_1+u_2 ,v_1+v_2)
=\frac{ \vartheta\left[\begin{array}{cc} 0 \ 1 \\ 0 \ 1 \\ \end{array}\right](u_1+u_2 ,v_1+v_2)} 
{\vartheta\left[\begin{array}{cc} 0 \ 0 \\ 0 \ 0 \\ \end{array}\right](u_1+u_2, v_1+v_2)} 
=\frac{G_1[\begin{array}{cc} 0 \ 1 \\ 0 \ 0 \\ \end{array}]
 G[\begin{array}{cc} 0 \ 1 \\ 0 \ 1 \\ \end{array}] }
{G[\begin{array}{cc} 0 \ 1 \\ 0 \ 0 \\ \end{array}]
G[\begin{array}{cc} 0 \ 0 \\ 0 \ 0 \\ \end{array}]}    ,
\label{A5}
\end{eqnarray}
\begin{eqnarray}
&&\hskip-15mm F[\begin{array}{cc} 0 \ 1 \\ 1 \ 0 \\ \end{array}](u_1+u_2 ,v_1+v_2)
=\frac{ \vartheta\left[\begin{array}{cc} 0 \ 1 \\ 1 \ 0 \\ \end{array}\right](u_1+u_2 ,v_1+v_2)} 
{\vartheta\left[\begin{array}{cc} 0 \ 0 \\ 0 \ 0 \\ \end{array}\right](u_1+u_2, v_1+v_2)}
=\frac{G_1[\begin{array}{cc} 0 \ 1 \\ 0 \ 0 \\ \end{array}]
 G[\begin{array}{cc} 0 \ 1 \\ 1 \ 0 \\ \end{array}] }
{G[\begin{array}{cc} 0 \ 1 \\ 0 \ 0 \\ \end{array}]
G[\begin{array}{cc} 0 \ 0 \\ 0 \ 0 \\ \end{array}]}   ,
\label{A6}
\end{eqnarray}
\begin{eqnarray}
&&\hskip-15mm F[\begin{array}{cc} 0 \ 1 \\ 1 \ 1 \\ \end{array}](u_1+u_2 ,v_1+v_2)
=\frac{ \vartheta\left[\begin{array}{cc} 0 \ 1 \\ 1 \ 1 \\ \end{array}\right](u_1+u_2 ,v_1+v_2)} 
{\vartheta\left[\begin{array}{cc} 0 \ 0 \\ 0 \ 0 \\ \end{array}\right](u_1+u_2, v_1+v_2)}
=\frac{G_1[\begin{array}{cc} 0 \ 1 \\ 0 \ 0 \\ \end{array}]
 G[\begin{array}{cc} 0 \ 1 \\ 1 \ 1 \\ \end{array}] }
{G[\begin{array}{cc} 0 \ 1 \\ 0 \ 0 \\ \end{array}]
G[\begin{array}{cc} 0 \ 0 \\ 0 \ 0 \\ \end{array}]}    ,
\label{A7}
\end{eqnarray}
\begin{eqnarray}
&&\hskip-15mm F[\begin{array}{cc} 1 \ 0 \\ 0 \ 0 \\ \end{array}](u_1+u_2 ,v_1+v_2)
=\frac{ \vartheta\left[\begin{array}{cc} 1 \ 0 \\ 0 \ 0 \\ \end{array}\right](u_1+u_2 ,v_1+v_2)} 
{\vartheta\left[\begin{array}{cc} 0 \ 0 \\ 0 \ 0 \\ \end{array}\right](u_1+u_2, v_1+v_2)} 
=\frac{ G_1[\begin{array}{cc} 1 \ 0 \\ 0 \ 0 \\ \end{array}] }
{G[\begin{array}{cc} 0 \ 0 \\ 0 \ 0 \\ \end{array}]}   ,
\label{A8}
\end{eqnarray}
\begin{eqnarray}
&&\hskip-15mm F[\begin{array}{cc} 1 \ 0 \\ 0 \ 1 \\ \end{array}](u_1+u_2 ,v_1+v_2)
=\frac{ \vartheta\left[\begin{array}{cc} 1 \ 0 \\ 0 \ 1 \\ \end{array}\right](u_1+u_2 ,v_1+v_2)} 
{\vartheta\left[\begin{array}{cc} 0 \ 0 \\ 0 \ 0 \\ \end{array}\right](u_1+u_2, v_1+v_2)} 
=\frac{G_1[\begin{array}{cc} 1 \ 0 \\ 0 \ 0 \\ \end{array}]
 G[\begin{array}{cc} 1 \ 0 \\ 0 \ 1 \\ \end{array}] }
{G[\begin{array}{cc} 1 \ 0 \\ 0 \ 0 \\ \end{array}]
G[\begin{array}{cc} 0 \ 0 \\ 0 \ 0 \\ \end{array}]}    ,
\label{A9}
\end{eqnarray}
\begin{eqnarray}
&&\hskip-15mm F[\begin{array}{cc} 1 \ 0 \\ 1 \ 0 \\ \end{array}](u_1+u_2 ,v_1+v_2)
=\frac{ \vartheta\left[\begin{array}{cc} 1 \ 0 \\ 1 \ 0 \\ \end{array}\right](u_1+u_2 ,v_1+v_2)} 
{\vartheta\left[\begin{array}{cc} 0 \ 0 \\ 0 \ 0 \\ \end{array}\right](u_1+u_2, v_1+v_2)}
=\frac{G_1[\begin{array}{cc} 1 \ 0 \\ 0 \ 0 \\ \end{array}]
 G[\begin{array}{cc} 1 \ 0 \\ 1 \ 0 \\ \end{array}] }
{G[\begin{array}{cc} 1 \ 0 \\ 0 \ 0 \\ \end{array}]
G[\begin{array}{cc} 0 \ 0 \\ 0 \ 0 \\ \end{array}]}   ,
\label{A10}
\end{eqnarray}
\begin{eqnarray}
&&\hskip-15mm F[\begin{array}{cc} 1 \ 0 \\ 1 \ 1 \\ \end{array}](u_1+u_2 ,v_1+v_2)
=\frac{ \vartheta\left[\begin{array}{cc} 1 \ 0 \\ 1 \ 1 \\ \end{array}\right](u_1+u_2 ,v_1+v_2)} 
{\vartheta\left[\begin{array}{cc} 0 \ 0 \\ 0 \ 0 \\ \end{array}\right](u_1+u_2, v_1+v_2)} 
=\frac{G_1[\begin{array}{cc} 1 \ 0 \\ 0 \ 0 \\ \end{array}]
 G[\begin{array}{cc} 1 \ 0 \\ 1 \ 1 \\ \end{array}] }
{G[\begin{array}{cc} 1 \ 0 \\ 0 \ 0 \\ \end{array}]
G[\begin{array}{cc} 0 \ 0 \\ 0 \ 0 \\ \end{array}]}     ,
\label{A11}
\end{eqnarray}
\begin{eqnarray}
&&\hskip-15mm F[\begin{array}{cc} 1 \ 1 \\ 0 \ 0 \\ \end{array}](u_1+u_2 ,v_1+v_2)
=\frac{ \vartheta\left[\begin{array}{cc} 1 \ 1 \\ 0 \ 0 \\ \end{array}\right](u_1+u_2 ,v_1+v_2)} 
{\vartheta\left[\begin{array}{cc} 0 \ 0 \\ 0 \ 0 \\ \end{array}\right](u_1+u_2, v_1+v_2)} 
=\frac{ G_1[\begin{array}{cc} 1 \ 1 \\ 0 \ 0 \\ \end{array}] }
{G[\begin{array}{cc} 0 \ 0 \\ 0 \ 0 \\ \end{array}]}    ,
\label{A12}
\end{eqnarray}
\begin{eqnarray}
&&\hskip-15mm F[\begin{array}{cc} 1 \ 1 \\ 0 \ 1 \\ \end{array}](u_1+u_2 ,v_1+v_2)
=\frac{ \vartheta\left[\begin{array}{cc} 1 \ 1 \\ 0 \ 1 \\ \end{array}\right](u_1+u_2 ,v_1+v_2)} 
{\vartheta\left[\begin{array}{cc} 0 \ 0 \\ 0 \ 0 \\ \end{array}\right](u_1+u_2, v_1+v_2)} 
=\frac{G_1[\begin{array}{cc} 1 \ 1 \\ 0 \ 0 \\ \end{array}]
 G[\begin{array}{cc} 1 \ 1 \\ 0 \ 1 \\ \end{array}] }
{G[\begin{array}{cc} 1 \ 1 \\ 0 \ 0 \\ \end{array}]
G[\begin{array}{cc} 0 \ 0 \\ 0 \ 0 \\ \end{array}]}    ,
\label{A13}
\end{eqnarray}
\begin{eqnarray}
&&\hskip-15mm F[\begin{array}{cc} 1 \ 1 \\ 1 \ 0 \\ \end{array}](u_1+u_2 ,v_1+v_2)
=\frac{ \vartheta\left[\begin{array}{cc} 1 \ 1 \\ 1 \ 0 \\ \end{array}\right](u_1+u_2 ,v_1+v_2)} 
{\vartheta\left[\begin{array}{cc} 0 \ 0 \\ 0 \ 0 \\ \end{array}\right](u_1+u_2, v_1+v_2)} 
=\frac{G_1[\begin{array}{cc} 1 \ 1 \\ 0 \ 0 \\ \end{array}]
 G[\begin{array}{cc} 1 \ 1 \\ 1 \ 0 \\ \end{array}] }
{G[\begin{array}{cc} 1 \ 1 \\ 0 \ 0 \\ \end{array}]
G[\begin{array}{cc} 0 \ 0 \\ 0 \ 0 \\ \end{array}]}    ,
\label{A14}
\end{eqnarray}
\begin{eqnarray}
&&\hskip-15mm F[\begin{array}{cc} 1 \ 1 \\ 1 \ 1 \\ \end{array}](u_1+u_2 ,v_1+v_2)
=\frac{ \vartheta\left[\begin{array}{cc} 1 \ 1 \\ 1 \ 1 \\ \end{array}\right](u_1+u_2 ,v_1+v_2)} 
{\vartheta\left[\begin{array}{cc} 0 \ 0 \\ 0 \ 0 \\ \end{array}\right](u_1+u_2, v_1+v_2)} 
=\frac{G_1[\begin{array}{cc} 1 \ 1 \\ 0 \ 0 \\ \end{array}]
 G[\begin{array}{cc} 1 \ 1 \\ 1 \ 1 \\ \end{array}] }
{G[\begin{array}{cc} 1 \ 1 \\ 0 \ 0 \\ \end{array}]
G[\begin{array}{cc} 0 \ 0 \\ 0 \ 0 \\ \end{array}]}    .
\label{A15}
\end{eqnarray}
In the right-hand side, we can express with hyperelliptic function just by
replacing\\
$\displaystyle{\vartheta\left[\begin{array}{cc} a \ c \\ b \ d \\ \end{array}\right](u_1, v_1)
\rightarrow F[\begin{array}{cc} a \ c \\ b \ d  \\ \end{array}](u_1 ,v_1)}$, 
$\displaystyle{\vartheta\left[\begin{array}{cc}a \ c \\ b \ d  \\ \end{array}\right](u_2, v_2)
\rightarrow F[\begin{array}{cc} a \ c \\ b \ d  \\ \end{array}](u_2 ,v_2)}$ and\\
$\displaystyle{\vartheta\left[\begin{array}{cc} a \ c \\ b \ d  \\ \end{array}\right](0, 0)
\rightarrow F[\begin{array}{cc} a \ c \\ b \ d  \\ \end{array}](0 ,0)}$, 
because $\displaystyle{\vartheta\left[\begin{array}{cc} * \ * \\ * \ * \\ \end{array}\right](u_1, v_1)}$
and $\displaystyle{\vartheta\left[\begin{array}{cc} * \ * \\ * \ * \\ \end{array}\right](u_2, v_2)}$ 
appear in the second order homogeneous way for each factor, so that we can rescale 
with the common factor.\\
Addition formulae Eq.(\ref{A1}) $\sim$ Eq.(\ref{A15}) with Eq.(\ref{B1})-Eq.(\ref{B19})
in Appendix B are numerically checked by REDUCE.

\newpage
%
%
%
%
%
\section{Duplication formulae of hyperelliptic theta functions}
\setcounter{equation}{0}
\begin{eqnarray}
&&\hskip-15mm G[
\right](2u_2 , 2v_2)   ,
\nonumber\\
\label{B19}
\end{eqnarray}
As duplication formulae, Eq.(\ref{B1})-Eq.(\ref{B19}) are numerically checked by REDUCE.
\newpage
%
%
%
%
%
\section{\large Functional Relations}
\setcounter{equation}{0}
\begin{eqnarray}
&&\hspace{-20mm} \varTheta[%
\right](u ,v) \right)  .
\label{C28}
\end{eqnarray}
Functional relations Eq.(\ref{C1}) $\sim$ Eq.(\ref{C28}) are numerically checked by REDUCE.
%
%
%
%
%
%
%
%
\newpage
\setcounter{equation}{0}
\section{\large Constants Relations}
\begin{eqnarray}
&&\hskip-25mm 
\varTheta\left[%
\right](0, 0)   }\  .
\label{D16}
\end{eqnarray}
Constants relations Eq.(\ref{D1}) $\sim$ Eq.(\ref{D16}) are numerically checked by REDUCE.
%
%
%
%
%
%
%
\newpage
\setcounter{equation}{0}
\section{\large $SO(3)$ Structure of Elliptic Addition Formulae}
Lie group is closed by the multiplication of elements as  
addition formulae of Lie group.
While algebraic functions such as trigonometric/hyperbolic/elliptic
functions are closed by addition formulae of algebraic functions. 
Then we expect that there is the relation
between Lie group and algebraic functions through addition 
formulae.
For the genus one case, this can be realized by representing the addition formula 
of the spherical trigonometry (SO(3) Lie group structure) by the elliptic 
function.\\ 
Addition formulae of the spherical trigonometry can be rewritten  
of the form
%
%
 \begin{eqnarray}
  &&\hspace{-20mm} 
  \left( \begin{array}{ccc}
  1 & 0 & 0 \\
  0 & \cos(A_3) & \sin(A_3) \\
  0 & -\sin(A_3) & \cos(A_3)
 \end{array} \right)
  \left( \begin{array}{ccc}
  \cos(a_2) & \sin(a_2) & 0 \\
  -\sin(a_2) & \cos(a_2) & 0 \\
  0 & 0 & 1
 \end{array} \right)
 \left( \begin{array}{ccc}
  1 & 0 & 0 \\
  0 & \cos(A_1) & \sin(A_1) \\
  0 & -\sin(A_1) & \cos(A_1)
 \end{array} \right)   \nonumber\\
 &&\hspace{-20mm} 
=\left( \begin{array}{ccc}
  \cos(a_1) & \sin(a_1) & 0 \\
  -\sn(a_1) & \cos(a_1) & 0 \\
  0 & 0 & 1
 \end{array} \right )
\left( \begin{array}{ccc}
  1 & 0 & 0 \\
  0 & \cos(\pi-A_2) & \sin(\pi-A_2) \\
  0 & -\sin(\pi-A_2) & \cos(\pi-A_2)
 \end{array} \right) 
 \left( \begin{array}{ccc}
  \cos(a_3) & \sin(a_3) & 0 \\
  -\sin(a_3) & \cos(a_3) & 0 \\
  0 & 0 & 1
 \end{array} \right )  .
\nonumber
\end{eqnarray}
The right-hand side is the Euler angle parametrization for the sphere and the left-hand side 
is the dual Euler angle parametrization, and the Euler angle parametrization and 
the dual Euler angle parametrization becomes equal.
\begin{figure}[htb]
 \centering
  \includegraphics{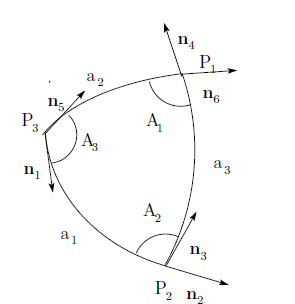}
  \caption{Euler angle parametrization and dual Euler angle parametrization}
\end{figure}%
\\ 
The relation of the $11$-component of  (the left-hand side)-(right-hand side) becomes 
$$\cos(a_2)-\cos(a_1) \cos(a_3)-\cos(A_2) \sin(a_1) \sin(a_3)=0, $$
and the relation of the $33$-component of  (the left-hand side)-(right-hand side) becomes 
$$\cos(A_2)+\cos(A_1) \cos(A_3)-\cos(a_2) \sin(A_1) \sin(A_3)=0, $$ 
which are some of spherical $\cos$ relations.
The relation of the $13$-component of  (the left-hand side)-(right-hand side) becomes 
$$-\sin(a_1) \sin(A_2)+\sin(a_2) \sin(A_1)=0, $$ 
which is one of spherical $\sin$ relations.
The relation of the $12$-component of  (the left-hand side)-(right-hand side) becomes 
$$-\cos(a_1) \sin(a_3)+\cos(A_1) \sin(a_2)+\cos(A_2)\sin(a_1) \cos(a_3)=0, $$
and the relation of the $23$-component of  (the left-hand side)-(right-hand side) becomes 
$$-\cos(a_1) \sin(A_2)+\cos(A_1) \sin(A_3)+\cos(A_3) \sin(A_1)\cos(a_2)=0, $$
and the relation of the $22$-component of  (the left-hand side)-(right-hand side) becomes 
$$\sin(a_1) \sin(a_3)-\sin(A_1) \sin(A_3)+\cos(A_2)\cos(a_1) \cos(a_3)
+\cos(a_2) \cos(A_1) \cos(A_3)=0,$$ 
those are shown to satisfy by using spherical trigonometric relations. 
Above relations are parameterized by elliptic functions of the 
form\cite{Shigemoto1,Shigemoto2} 
\begin{eqnarray}
  &&\hspace{-10mm} 
  \left( \begin{array}{ccc}
  1 & 0 & 0 \\
  0 & \dn(u_3) & k \sn(u_3) \\
  0 & -k \sn(u_3) & \dn(u_3)
 \end{array} \right)
  \left( \begin{array}{ccc}
  \cn(u_2) & \sn(u_2) & 0 \\
  -\sn(u_2) & \cn(u_2) & 0 \\
  0 & 0 & 1
 \end{array} \right)
 \left( \begin{array}{ccc}
  1 & 0 & 0 \\
  0 & \dn(u_1) & k \sn(u_1) \\
  0 & -k \sn(u_1) & \dn(u_1)
 \end{array} \right)    \nonumber\\
 &&\hspace{-10mm} 
 =\left( \begin{array}{ccc}
  \cn(u_1) & \sn(u_1) & 0 \\
  -\sn(u_1) & \cn(u_1) & 0 \\
  0 & 0 & 1
 \end{array} \right )
\left( \begin{array}{ccc}
  1 & 0 & 0 \\
  0 & \dn(u_2) & k \sn(u_2) \\
  0 & -k \sn(u_2) & \dn(u_2)
 \end{array} \right) 
 \left( \begin{array}{ccc}
  \cn(u_3) & \sn(u_3) & 0 \\
  -\sn(u_3) & \cn(u_3) & 0 \\
  0 & 0 & 1
 \end{array} \right )   , 
\nonumber 
\end{eqnarray}
with $u_2=u_1+u_3$. This is satisfied by using addition formulae of elliptic functions.
This is the Yang-Baxter equation of the Ising model\cite{Shigemoto1,Shigemoto2},
that is, the $SO(3) \cong SU(2)/Z_2$  integrability condition parametrized 
by elliptic functions.\\
%
%
%
%
%
%
%
\end{document}